\documentclass[12pt]{article}
\usepackage{amssymb,amsmath}
\usepackage{graphicx}
\usepackage[perpage,symbol]{footmisc}
\def\ba#1{\begin{array}{#1}}
\def\ea{\end{array}}
\def\beq#1{\begin{equation}\label{#1}}
\def\eeq{\end{equation}}
\renewcommand{\Re}{{\rm Re}\,}
\renewcommand{\Im}{{\rm Im}\,}
\newcommand{\RR}{\mathbb{R}}
\newcommand{\Rplus}{\RR_+}
\newcommand{\NN}{\mathbb{N}}
\newcommand{\CC}{\mathbb{C}}
\newcommand{\ZZ}{\mathbb{Z}}
\newcommand{\HH}{\mathbb{H}} 
\newcommand{\Can}{\mathcal{C}} 
\newcommand{\PCan}{\Can\times\Can} 

\newcommand{\cP}{\mathcal{P}} 
\newcommand{\Cp}{\CC_+} 
\newcommand{\Ray}[1]{\Lambda_{#1}}

\newcommand{\unl}{\underline}
\newcommand{\intl}{\int\limits}
\def\Lsp#1{L^{#1}} 
\newcommand{\Lp}{\Lsp{p}}
\newcommand{\Lpp}{\Lsp{p'}}

\newcommand{\f}{\frac}
\def\num#1{{\rm(\ref{#1})}}

\newcommand{\lam}{\lambda}
\newcommand{\al}{\alpha}
\newcommand{\be}{\beta}
\newcommand{\om}{\omega}
\newcommand{\eps}{\varepsilon}
\newcommand{\ph}{\varphi}
\newcommand{\const}{\mathrm{const}\,}

\newcommand{\bom}{\overline{\om}}
\newcommand{\tlf}{\tilde f}
\newcommand{\ga}{\gamma}
\newcommand{\Ga}{\Gamma}

\newcommand{\Ds}{E} 
\newcommand{\Df}{d} 
\newcommand{\Mes}{\mathfrak{M}} 
\newcommand{\Lt}{\mathcal{L}} 
\newcommand{\St}{\mathcal{S}} 
\newcommand{\Ht}{\mathcal{H}} 
\newcommand{\Ft}{\mathcal{F}} 

\newcommand{\Lg}{\Lt_\ga}
\newcommand{\Lr}{\Lt_\theta}

\newcommand{\Sta}{\St_\theta}
\newcommand{\Mon}{\mathrm{Mon}} 

\newcommand{\MonCP}{\Mon(\Cp)}
\newcommand{\ConvCP}{\mathrm{Conv}(\Cp)} 
\newcommand{\Wp}{\mathfrak{W}} 
\newcommand{\KPcurve}{K_3}
\newcommand{\KScurve}{K_4}
\newcommand{\KLcurve}{K_5}
\newcommand{\Kmon}{K_6}
\newcommand{\Kconv}{K_7}
\newcommand{\Kgraph}{K_8}
\newcommand{\Kgraphy}{\hat{K}_8}
\newcommand{\Kbox}{K_{9}}
\newcommand{\KmaxPW}{K_{10}}
\newcommand{\Kmaxang}{K_{11}}
\newcommand{\Kvcomb}{K_{12}}

\def\dst{\displaystyle}
\def\eop{\hfill$\Box$}
\def\rem#1{}

\newcounter{ex} 
\newcommand{\ex}[1]{\refstepcounter{ex}%
\bigskip\noindent{\bf Example \arabic{ex}}. \label{#1}}

\newtheorem{defn}{Definition}
\newtheorem{thm}{Theorem}
\newtheorem{lem}{Lemma}

\newcommand{\pf}{\noindent{\it Proof}. }

\newcommand{\newsectnum}{\setcounter{equation}{0}  \setcounter{lem}{0}}

\begin{document}

\begin{center}
{\Large 
Hausdorff-Young type theorems for the Laplace transform restricted to
a ray or to a curve in the complex plane
}
\\[4ex]{\large
{ Anatoli Merzon}\footnote{ 
Instituto de F\'isica y Matem\'aticas,
Universidad Michoacana de San Nicol\'as de Hidalgo,
Morelia, Michoac\'an, M\'exico.
E-mail:	
{\tt anatoli@ifm.umich.mx}}
and
{Sergey Sadov}
\footnote{
Department of Mathematics and Statistics, Memorial University of Newfoundland,
St.~John's, Canada.
E-mail:	
{\tt sergey@mun.ca}}
}
\end{center}

%
%
%

%
%
%


\begin{quote} 
\vspace{5ex}
\centerline{\bf Abstract}

\noindent 
Let $\mathcal{L} f$ 
be the Laplace transform of a function $f\in X = L^p(0,\infty)$. For some classes of Borel measures $\mu$ in the half-plane $\Re z>0$ we prove the Hausdorff-Young property $\|\mathcal{L} f\|_Y \leq K\|f\|_X$,
where $Y = L^{p'}(d\mu)$, $p$ and
$p'$ are conjugate exponents, $1\leq p\leq 2$, and the
constant $K$ is uniform over the class. Particular attention is paid to the case $d\mu=ds_{\gamma}$, the arclength measure on a rectifiable curve $\gamma$ belonging to a suitable family.
Examples are the family of convex curves and the family of curves consisting of convex components enclosed in a horizontal row of boxes with bounded height-to-width ratio. In addition, for rays $\mathrm{arg}\, z = \theta$, 
the exact $L^2$ norm of the operator $f\mapsto \mathcal{L}f(re^{i\theta})$ is given and an analog of the
Hausdorff-Young inequality in Lorentz spaces is obtained in the case of ``wrong'' exponents $p > 2$.

\bigskip\bigskip
\hspace*{-28pt}{\bf Keywords}: 
Laplace transform, Fourier transform, Hausdorff-Young inequality, Paley-Wiener theory, arclength measure, Lorentz spaces.
\end{quote}

\newpage


\tableofcontents
\contentsline{section}{Figures}{61}

\newpage
\section{Introduction}

The Laplace transform of a function $f(t)$ defined on $\RR_+=(0,\infty)$
(from an appropriate class)
is an analytic function in the right half-plane $\Cp=\{z=x+iy 
|\;x>0\}$,
\beq{Lt}
 \Lt f(z)=\int_0^\infty f(t)e^{-zt}\,dt.
\eeq
Let $\Ray{\theta}$ be the ray $z=\rho e^{i\theta}$, $\rho>0$,
$\theta=\const\in[-\pi/2,\pi/2]$. The Laplace transform observed along the ray $\Ray{\theta}$ is
\beq{Lr}
 \Lr f(\rho)=\int_0^\infty f(t)e^{-\rho t e^{i\theta}}\,dt.
\eeq

Generalizing, suppose $\ga\subset\Cp$ is a rectifiable curve parametrized by arclength $z=z(s)$, $s\in I\subset \RR$.
We consider the Laplace transform observed along $\gamma$ as a function of the parameter $s$:
\beq{Lg}
 \Lg f(s)=\int_0^\infty f(t)e^{-z(s)t}\,dt,\qquad s\in I.
\eeq

The subject of this paper, which extends our earlier publication \cite{SM-CMA10},
is
$L^p$ estimates for the operators $\Lr$ and operators $\Lg$ corresponding to various families of curves. 
It turns convenient to treat arclength measures $|dz(s)|$ on many suitable curves as examples
of what we call well-projected measures on $\Cp$. A part
of this study is devoted to such measures and to  corresponding weighted $L^p$ estimates for the Poisson and Cauchy integrals.

\smallskip
Let us detail the plan of the paper and sketch the results.

\smallskip
In Sect.~\ref{sec:ray} we show that the operator $\Lr$ is 
a bounded operator from $L^p(\RR_+)$ to $\Lpp(\RR_+)$ for
$1\leq p\leq 2$, uniformly in $\theta$. Here $p$ and $p'$ are conjugate exponents:
\beq{lexp}
 p^{-1}+p'^{-1}=1, \qquad 1\leq p\leq p'\leq \infty.
\eeq
This relation between $p$ and $p'$ is always assumed in the sequel.

Theorem \ref{lray2} determines the exact norm of the operator $\Lr$ in $L^2(\RR_+)$.
The proof given in Sect.~\ref{prray} is based on the explicit spectral decomposition of
the integral operator $\Lr^*\,\Lr$ with Hermite-symmetric kernel
$
k(x,y)=(xe^{-i\theta}+y e^{i\theta})^{-1}
$.

Theorem \ref{lrayp} asserts $\Lp(\RR_+)\to \Lpp(\RR_+)$ boundedness, uniformly in $\theta$, of $\Lr$ for $1\leq p\leq 2$.
It easily follows from Theorem \ref{lray2} by the Riesz-Thorin interpolation theorem.
The closer $|\theta|$ is to $\pi/2$ in (\ref{Lr}), the greater role in the estimates is due to oscillations
of the exponential factor $\exp(-\rho t e^{i\theta})$, as opposed to its decay in magnitude.
The limiting cases $\theta=\pm\pi/2$ correspond to the classical Hausdorff-Young theorem for the Fourier transform.

In Sect.~\ref{lpneg} we discuss the operators $\Lr$  in the spaces $\Lp$ with ``wrong'' exponents $p>2$.
For such values of $p$, neither the Fourier transform (or its half-line version $\Lt_{\pm\pi/2}$)
nor the Laplace transform $\Lt_{0}$ are bounded maps $L^p\to L^{p'}$. This fact for the Laplace transform,
while long known, is surprisingly little mentioned in the literature. We have included a counterexample  
(Lemma~\ref{lraywneg}). In contrast, Theorem~\ref{lrayw} states a positive result: $\Lr$ is a bounded map between the Lorentz spaces $L^{p,r}$ and $L^{p',r}$ if $|\theta|<\pi/2$.


Sections~\ref{sec:curve} and \ref{sec:app} are full of technicalities, but their objective is clear:
Hausdorff-Young type theorems for operators $\Lg$ with norm estimates uniform over different classes of curves in $\Cp$.
Let us review contents of these sections in 
motivational rather than in sequential order.

The target of Section~\ref{sec:curve} in short is Theorem~\ref{lmasterthm},
a ``Master Theorem''.
Theorems in Section~\ref{sec:app} are its corollaries, with exception of the last Theorem~\ref{thm:vertcomb}. 

Among these corollaries, the first and easiest to understand 
is Theorem~\ref{thm:unifcurv}.
It covers in particular the class of convex curves
(boundaries of convex domains in $\Cp$).
In the special case of hyperbolas
\beq{hyp}
\f{x^2}{a^2}-\f{y^2}{b^2}=1,\qquad x>0,
\eeq
the $L^2$ estimates for $\Lg$ have been proven in our paper
\cite{SM-CMA10}. A motivation came from mathematical physics
--- a study of the Neumann problem for the Helmholtz operator in a plane angle,  
%
%
see \cite[Lemma 7.2]{ZM} 
and \cite{Merzon09}, ---
and prompted a more general inquiry into the Laplace transform observed along curves.

In spite of a much more narrow focus of \cite{SM-CMA10},
all ingredients of the proof of the master Theorem~\ref{lmasterthm} have been already exposed there. 
In the proof of Lemma~3.2 of \cite{SM-CMA10} the key role 
was played by a quasi-triangle inequality, with length of
a subset of a plane curve on one side and lengths of the $x$- and $y$-projections of that subset on the other. Axiomatizing, we call curves that possess such a property well-projected curves. Abstracting further,
we come to the notion of a well-projected measure, which
generalizes the arclength on a well-projected curve.

The only place, but a critical one, where well-projectedness
enters into proof, is the end of proof of part (a) of 
Theorem~\ref{pcurve} (p.~\pageref{pg:eoppcurve}).

Theorems~\ref{pcurve} and \ref{ccurve}
of Sect.~\ref{sec:cauchy}, 
on which our proof of Theorem~\ref{lmasterthm} is based, 
assert $L^p$ estimates with well-projected weights for the Poisson integral and the Cauchy integral, respectively.
They may be of interest in their own right. 

As mentioned earlier, most results in Sect.~\ref{sec:app}
are corollaries of Theorem~\ref{lmasterthm}. However, before referring to Theorem~\ref{lmasterthm} one must ensure that the measures involved are well-projected. Sometimes it is quite obvious and sometimes requires a certain effort.
That part of the work --- checking well-projectedness and
determining the corresponding constants --- is mostly done
in Sect.~\ref{sec:wpexamples}. Thus, Theorem~\ref{thm:unifcurv} refers to Examples \ref{ex:moncurve} and \ref{ex:convcurve}, Theorem~\ref{thm:lipcurve} --- to Example~\ref{ex:lipcurve},
Theorem~\ref{thm:boxcurve} --- to Example~\ref{ex:boxedcurve},
and Theorem~\ref{restrangmax} --- to Example~\ref{ex:maxang}.

A ``Maximal Paley-Wiener Theorem''~\ref{lmaxPW}, not yet mentioned, stands out in that it deals with a nonlinear maximal operator $\Lt^*$ instead of
linear operators $\Lg$ in other theorems of Section~\ref{sec:app}. It is similar to the Hardy-Littlewood theorem for the maximal radial function
of a function in the Hardy space $H^p$ on the unit disk \cite[Th.~1.9]{Duren1970}.  
Whereas the Paley-Wiener theory asserts that
\beq{pw-supint}
 \sup\limits_{x>0} \;\int\limits_{-\infty}^\infty |\Lt f(iy)|^{p'}\,dy \,\leq \, C \|f\|_p^p,
\eeq
in Theorem \ref{lmaxPW} we first take the supremum
and then integrate:
\beq{pw-intsup}
 \int_{-\infty}^\infty \left(\sup_{x>0} 
|\Lt f(x+iy)|\right)^{p'}\,dy \,\leq \, \tilde C \|f\|_p^p.
\eeq
By a simple trick the estimate (\ref{pw-intsup}) is equivalent to an uniform estimate for $\|\Lg\|_{L^{p}\to L^{p'}}$ over a family of discontinuous curves consisting
of vertical segments with non-overlapping $y$-projections.
    
Theorem~\ref{restrangmax} is of a similar nature.
However, it is not a 
kind of
result in a nice, ``final'' form. Rather we perceive it as a prototype (with extra conditions which can likely be dropped or at least significally weakened) of a ``Maximal Angular Inequality'', see (\ref{maxangpure}),
which would be in the same relation to the uniform in $\theta$ estimate for $\|\Lr\|_{L^p\to L^{p'}}$ from
Theorem~\ref{lrayp} as the inequality~(\ref{pw-intsup}) is
to (\ref{pw-supint}).
Theorem~\ref{lmasterthm} seems to be insufficient to
obtain the Maximal Angular Inequality in a desired generality.

The last Theorem~\ref{thm:vertcomb}
demonstrates one example where Theorem~\ref{lmasterthm}
is insufficient, too, but an enhancement of technique
used in its proof yields the result.

\medskip


Tools employed in this work belong to the ``old'' harmonic analysis: 
Riesz-Thorin and Marcinkiewicz interpolation theorems, Hilbert transform, Hardy-Littlewood
maximal theorem, with Lorentz spaces and Calderon's extension of Marcinkiewicz's theorem being about the most
recent.
Most results 
could have been obtained 
in the 1920s-30s; however we were not able to locate even the fairly elementary Theorems \ref{lray2}
 and \ref{lrayp} (for $\theta\neq 0,\pm\pi/2$) in the literature.

\subsection*{Notation and conventions}

\medskip
\noindent{\em Constants}\\
The letter $C$ with or without subscripts denotes a positive constant in inequalities, always within a local scope (such as proof of a lemma or a theorem). The dependence on parameter(s) 
may be omitted, but is
indicated explicitly where necessary. 

The subscripted letter $K$ denotes the best constants (norms) in statements of theorems. 
The scope of each $K_i$ is global and the $K_i$'s may depend on parameter(s). For example,
$K_2(\theta,p)$, defined in Theorem~2, is the same throughout the paper.   
Precise values of the $K_i$'s are in most cases unknown to the authors; available bounds are given.

\medskip
\noindent{\em Measure}\\
%
The Lebesgue measure on the real line or on its subsets is denoted $dt$, $dx$, etc. The arclength measure 
on rectifiable curves in the complex plane $\CC$ is 
denoted $ds$ or $|dz|$.

If $E$ is a subset of $\RR$ or of a rectifiable curve in $\CC$, then its measure in the applicable sense as above (i.e.\ length) is denoted $|E|$.
When dealing with more general measures, we use notation
$\mu(E)$; here $E$ is a Borel subset of $\CC$.

%
For $f$ a real- or complex-valued function and $\lambda>0$ 
we write
\beq{Ds}
\Ds_f(\lambda)=\{x:\; |f(x)|>\lambda\}.
\eeq
It is a subset of the domain of $f$.
The {\em distribution function}\ of $f$ is
\beq{Df}
 \Df_f(\lambda)=|E_f(\lambda)|
\qquad \mbox{\rm or}\qquad
\Df_f(\lambda)=\mu(E_f(\lambda)),
\eeq
depending on a relevant measure, which should be
clear from the context. 

\medskip
\noindent{\em Norms}\\
For a measurable function whose domain is an interval $I$ of the real line (resp.\ a rectifiable curve $\ga$ in $\CC$ or a more general set equipped with measure $\mu$),
the $L^p$-norm with respect to the measure $dt$ (resp.\ $ds$
or $\mu$) is denoted $\|f\|_{\Lp(I)}$ (resp.\ $\|f\|_{\Lp(\ga)}$ or $\|f\|_{\Lp(d\mu)})$,
or simply $\|f\|_p$ if the context ensures clarity.  

\smallskip
The weak space $L^{p,\infty}(I)$ ($1\leq p<\infty$)
consists of measurable functions $f$ on $I$ for which there exists
$C>0$ such that
\beq{weakLp}
\Df_f(\lambda)\leq (C/\lambda)^{p},\quad \forall\lambda>0,
\eeq
and the weak $p$-norm of $f$ is $\|f\|_{L^{p,\infty}(I)}=\inf C$. We write $\|f\|_{p,\infty}$ if $I$ is implicit.
(In fact, $\|\cdot\|_{p,\infty}$ is not a norm but a quasinorm, see e.g.\ \cite[Sect.~1.3]{BL}.)

\smallskip
More general Lorentz (quasi)norms $\|\cdot\|_{p,r}$
appear in Sect.~\ref{lpLorentz} and nowhere else in the
paper. The relevant notation is located there.

\smallskip
If an operator $A$ (not necessarily linear) carries functions from $\Lp(I)$ into $L^q(J)$,
then  
$$
 \|A\|_{\Lp(I)\to L^q(J)}=\inf\left\{C\,\left|\;\, \|Af\|_{L^q(J)}\leq C\|f\|_{L^p(I)},\quad\forall f\in\Lp(I)\right.\right\}.
$$
Similar notation is used for norms of operators whose 
target space is a weak space $L^{q,\infty}(J)$.
%
Whenever it is safe to do so, we abbreviate $\,\|A\|_{\Lp(I)\to L^q(J)}\,$ to $\,\|A\|_{p,q}\,$. 
If $A$ is an operator of weak type $(p,q)$, that is, $A:\;\Lp(I)\to L^{q,\infty}(J)$, then
we abbreviate $\,\|A\|_{\Lp(I)\to L^{q,\infty}(J)}\,$ to $\,\|A\|_{p,q*}\,$

The notation $\|f\|_{p,\infty}$ will stand for the
weak $p$-norm of a function $f$, while $\|\mathcal{A}\|_{p,\infty}$
will denote the norm of an operator $\mathcal{A}$ acting from $L^p$ to $L^\infty$. Despite a minor conflict of notation, the context should suffice to prevent confusion.


\section{Hausdorff-Young type theorems for 
the La\-place transform observed along a ray}

\label{sec:ray}
\newsectnum

\subsection{$\Lp\to\Lpp$ boundedness, $1\leq p\leq 2$}

Fix $\theta\in[-\pi/2,\pi/2]$ and consider
the operator $\Lr:\,f(t)\mapsto\Lr f(\rho)$ defined by the formula (\ref{Lr}).
Note that $\Lt_0$ is the 
Laplace transform on the positive real half-line,
while $\Lt_{\pm\pi/2}$ are the restrictions onto a real half-line (positive or negative)
of the Fourier transform of the function $f$ with support in $[0,\infty)$.

It is easy to deduce $L^p\to L^{p'}$ boundedness of $\Lr$ for $\theta\in(-\pi/2,\pi/2)$ 
from the known $L^p\to L^{p'}$ boundedness of $\Lt_0$
\cite[Th.~352]{HLP}. Indeed, we have
\beq{Ltheta0}
 |\Lr f(\rho)|\leq \int_0^\infty e^{-t\rho \cos\theta}|f(t)|\,dt
= \Lt_0 |f|(\rho\cos\theta). 
\eeq     
Raising to the power $p'$ and integrating, we obtain
$$
\|\Lr f\|_{p'}\leq (\cos\theta)^{-1/p'}\,\|\Lt_0 |f|\;\|_{p'} 
$$
Therefore, 
\beq{Lpnormtheta0}
\|\Lr\|_{p,p'}\leq (\cos\theta)^{-1/p'}\,\|\Lt_0\|_{p,p'}.
\eeq
For estimates of this type in $L^p$ spaces with power weights
consult \cite[\S~4.2]{Sed2005}. 
Our goal is to obtain a stronger result: a norm estimate uniform in $\theta$. 

Note first that $\Lr$ is a bounded operator from
$L^1(\RR_+)$ to $L^\infty(\RR_+)$ with norm independent of $\theta$:
\beq{lray1}
\|\Lr\|_{1,\infty}=\sup_{t,\rho>0}|e^{-\rho t e^{i\theta}}|=1.
\eeq
Next we develop an $L^2$-theory of $\Lr$.   
Like in the $L^2$ theory of the Fourier transform, the first step is to establish
$L^2$-norm boundedness of $\Lr$ defined on $L^1\cap L^2$.
This step does not readily follow neither from Cauchy-Schwarz, nor from Schur's test, but it could be bypassed
by referring to (\ref{Lpnormtheta0}) with $p=2$. However, 
a more economical way to obtain 
a desired improvement over (\ref{Lpnormtheta0}) is to analyse the operator $\Lr|_{L^1\cap L^2}$ from scratch.
Then of course $\Lr$ is extended by continuity to a bounded operator in $L^2$ with an uniform in $\theta$ norm estimate.

\begin{thm}\label{lray2}
(a)  For any $\theta\in [-\pi/2,\pi/2]$ and any $f\in L^1(\RR_+)\cap L^2(\RR_+)$ the inequality
 \beq{Lr2}
  \|\Lr f\|_2\leq K_1(\theta)\|f\|_2
 \eeq
 holds. Thus $\Lr$ extends to the operator $L^2(\RR_+)\to L^2(\RR_+)$ by continuity.
\\[0.5ex]
(b)
 The best possible constants
\beq{K1theta}
 K_1(\theta)=\|\Lr\|_{2,2}
\eeq
in (\ref{Lr2}) are:
for $|\theta|<\pi/2$
 \beq{Ltnorm}
  K_1(\theta)=\sqrt{\pi}\,(1-\mu)^{({1-\mu})/{4}}\,(1+\mu)^{({1+\mu})/{4}},
  \qquad
 \mu=\f{|\theta|}{\pi/2},
 \eeq
 and for $|\theta|=\pi/2$
\beq{L2normpi}
K_1(\pm\pi/2)=\sqrt{2\pi}.
\eeq
(c)
 The even function $K_1(\theta)$ is increasing and continuous on $[0,\pi/2]$;
 in particular, the uniform norm estimate in $L^2$ for the operators $\Lr$ is
 \beq{L2norm}
 \sup_{|\theta|\leq \pi/2} K_1(\theta)
 =\sqrt{2\pi}.
 \eeq
\end{thm}

The proof in the case $|\theta|<\pi/2$ is based on the spectral analysis of the operators
$\Lr^*\Lr$ and will be given in the next section.
Here let us comment on the limiting cases $\theta=\pm\pi/2$, where
the operator $\Lr$ is the restriction of the Fourier transform
to the positive or negative real half-line.
The estimate
\beq{Ltriv}
\|\Lt_{\pm\pi/2}\|_{2,2}\leq\sqrt{2\pi}
\eeq
 follows from Plancherel's theorem.

We want to prove the equality (\ref{L2normpi}).
To be definite, let $\theta=\pi/2$.

Take a nonzero function $f(t)\in L^2(\RR)$  supported on $[0,\infty)$.
Consider the family of functions $f_\xi(t)=f(t)e^{-i\xi t}$, $\xi\in\RR$.
If $\Ft f(\rho)=\int_0^\infty e^{it\rho} f(t)\,dt\,$ is the Fourier transform of $f$, then
the Fourier transform of $f_\xi$ is $\Ft f_\xi(\rho)=\Ft f(\rho-\xi)$. Letting $\xi\to+\infty$,
we obtain
$$
\ba{c}\dst
 \|\Lt_{\pi/2} f\|^2_{L^2(\RR_+)}=\int_0^\infty |(\Ft f)(\rho-\xi)|^2\,d\rho=
\int_{-\xi}^\infty |\Ft f(\rho)|^2\,d\rho
\\[3ex]\dst
\to \|\Ft f\|_{L^2(\RR)}=\sqrt{2\pi}\|f\|_2.
\ea
$$
Thus the constant $\sqrt{2\pi}$ in (\ref{Ltriv}) is best possible.
\eop

\medskip
In the same way as the Plancherel theorem and the trivial $L^1\to L^\infty$ estimate
imply the Hausdorff-Young theorem by means of interpolation,
a $\theta$-uniform analog of the Hausdorff-Young now follows.

\begin{thm}\label{lrayp}
 Suppose $1\leq p\leq 2$ and $\theta\in [-\pi/2,\pi/2]$.
The operator $\Lr$ defined initially on $L^1(\RR_+)\cap L^p(\RR_+)$ extends to a bounded
 operator from $\Lp(\RR_+)$ to $\Lpp(\RR_+)$, where the relation (\ref{lexp}) is assumed, with norm
\beq{Lrp}
  \|\Lr\|_{p,p'}= K_2(\theta,p).
\eeq
Here
\beq{K2}
K_2(\theta,p)\leq K_1(\theta)^{2/p'}.
\eeq
Consequently, an uniform in $\theta$ bound is
\beq{Ltnormp}
  K_2(\theta,p)\leq (2\pi)^{1/p'},\quad\forall\theta\in[-\pi/2,\pi/2].
\eeq
\end{thm}

\pf For the operator $\Lr$ defined on
$L^1(\RR_+)\cap L^2(\RR_+)$ we have the norm estimates (\ref{lray1}) and 
(\ref{Lr2}).
By the Riesz-Thorin interpolation theorem 
\cite[Th.~1.1.1]{BL}
$\Lr$ is a bounded operator $\Lp(\RR_+) \to \Lpp(\RR_+)$ for any $p\in[1,2]$,  with norm
\beq{K2interp}
K_2(\theta,p)\leq 
\|\Lr\|_{1,\infty}^{2/p-1}\,\|\Lr\|_{2,2}^{2-2/p},
\eeq
which yields (\ref{K2}).
The uniform estimate (\ref{Ltnormp}) now follows from (\ref{L2norm}).
\eop

\medskip
\noindent
{\bf Remark}. Theorem~\ref{lrayp} gives only an upper
bound for $\|\Lr\|_{p,p'}$.

In the cases $\theta=0$ and $\theta=\pm\pi/2$ the simple estimate (\ref{K2}) that follows from the 
Riesz-Thorin theorem can be compared with known estimates derived by more refined methods. 

In the case $\theta=\pm\pi/2$ the Fourier transform) the precise constant for is known for all $p\in[1,2]$: 
\beq{BaBeLp}
K_2(\pm\pi/2,p)=(2\pi)^{1/p'}\,\left(\f{p^{1/p}}{p'^{1/p'}}\right)^{1/2}.
\eeq
This result is due to Beckner \cite{Beckner75} 
(and earlier to Babenko, for integer $p$).
In our situation the Fourier transform acts on functions with support in $\RR_+$, but the norm is not affected  
--- by the argument used above to prove that the constant in (\ref{L2norm}) is the same as in (\ref{Ltriv}).
The estimate (\ref{Ltnormp}) is inexact by about $16\%$: 
$$
\max\limits_{1\leq p\leq 2} \frac{(2\pi)^{1/p'}}{K_2(\pi/2,p)} \approx 1.158.
$$
(The maximum is attained at the smaller root $p\approx 1.192$ of the quadratic equation 
$p^2 - e^2p +e^2=0
$.)

For $|\theta|\neq\pi/2$ and $1<p<2$ we don't know precise values of $K_2(\theta,p)$. 
The estimate (\ref{K2}) says $K_2(0,p)\leq \pi^{1/p'}$. 
Hardy \cite[Th.~9]{Hardy33}, \cite[Th.~352]{HLP} showed that
\beq{HardyLp}
K_2(0,p)\leq(2\pi/p')^{1/p'}.
\eeq
Setterqvist \cite[Th.~2.2]{Sqv} substituted Beckner's sharp $L^p$ estimate for convolutions  into Hardy's proof and obtained
\beq{EricLp}
K_2(0,p)\leq \left(\pi(p-1)\right)^{1/p'}\,\left(p(2-p)\right)^{1/p-1/2}.
\eeq
The improvement over (\ref{K2}) given by (\ref{HardyLp}) and (\ref{EricLp}) is characterized as follows: 
$$
\max_{p} \frac{\pi^{1/p'}}{\rm r.h.s.\ of (\ref{HardyLp})} = e^{1/2e}\approx 1.202
$$
and
$$
\max_{p} \frac{\pi^{1/p'}}{\mbox{\rm r.h.s.\ of (\ref{EricLp}) }} \approx 1.355.
$$
(The maximum is attained at the positive root $p\approx 1.328$ of the quadratic equation 
$p^2 + (e-2)p -e=0
$.)

The known $L^2$ norm (\ref{K1theta}) of $\Lr$ 
and the method of \cite{Hardy33} yield the bound 
$$
K_2(\theta,p)\leq \left(\f{2\pi}{p'\cos\theta}\right)^{1/p'},
$$
which is better than (\ref{K2}) if $|\theta|<\arccos(2/p')$, but for a fixed $p\in (1,2]$ and $\theta\to\pi/2$ it diverges.
The same applies to a furhter enhancement (omitted) similar to  (\ref{EricLp}): it does nothing to the factor $(\cos\theta)^{-1/p'}$.
%
%
%

\subsection{Spectral decomposition in $L^2$: proof of Theorem \ref{lray2}, $|\theta|<\pi/2$
}
\label{prray}

\medskip
The proof of Theorem~\ref{lray2} consists of a series of lemmas.
Lemma \ref{l2absconv} and its corollary give boundedness of $\Lr$ in $L^2$ --- part (a) of the Theorem.
In Lemmas \ref{lls}--\ref{maxl} we study the self-adjoint operator
\beq{qf}
\Sta=\Lr^*\Lr.
\eeq
In Lemma~\ref{lemplancherel}, which is central to this section,
an explicit spectral decomposition of the operator $\Sta$ is constructed,
which readily yields $\|\Sta\|$ as the supremum of the spectrum of $\Sta$.
Thus we find
\beq{normLS}
\|\Lr\|_{2,2}=\|\Sta\|_{2,2}^{1/2},
\eeq
and obtain formula~(\ref{Ltnorm}).
The behaviour of the function $k(\theta)=\|\Sta\|^2_{2,2}$ is studied in Lemma~\ref{suplam}.
This furnishes a proof of the statement (c) of the Theorem.


\smallskip
In this section
$(\cdot,\cdot)$ denotes the scalar product in $L^2(\RR_+)$. The notation
$$
\omega=e^{i\theta}
$$
is used for brevity.

\medskip
In Lemma~\ref{l2absconv} 
and its 
proof the limits of integrations $0$ to $\infty$
are assumed throughout.

\medskip
We begin the proof
by showing that the bilinear form $(\Lr f,\Lr g)$
is  correctly defined and bounded on $L^2(\RR_+)\times L^2(\RR_+)$.

\begin{lem}
\label{l2absconv}
For any $f,g\in L^2(\RR_+)$ the inequality
\beq{tripleint}
\ba{ll}
\dst
|(\Lr f,\Lr g)|
&\dst
\leq
\int\kern-4pt\int\kern-4pt\int |f(y)\overline{g(x)}\, e^{-x\rho\bom-y\rho\om}|\,dx\,dy\,d\rho,
\\[3ex]
&\dst
\leq C(\theta)\|f\|_2\,\|g\|_2
\ea
\eeq
holds with constant
$C(\theta)$ independent of $f$ and $g$.
\end{lem}

\pf
By definition (\ref{Lr}) of the operator $\Lr$,
we have
\beq{ltripli}
(\Lr f,\Lr g)=
\int\kern-4pt\int\kern-4pt\int 
 f(y)\,\overline{g(x)}\, e^{-x\rho\bom-y\rho\om}\, dx\,dy\,d\rho.
\eeq
Therefore
$$
\ba{l}
|(\Lr f,\Lr g)|
\leq \dst
 \int\kern-4pt\int\kern-4pt\int 
 |f(y)|\,|g(x)|\, |e^{-x\rho\bom-y\rho\om}|\, dx\,dy\,d\rho
\\[3ex]
\qquad =\dst
 \int\kern-4pt\int\kern-4pt\int |f(y)|\,|g(x)|\, e^{-\rho(x+y)\cos\theta}\, dx\,dy\,d\rho
\\[3ex]
\qquad = \dst
 \int
 \kern-4pt
 \int
 \f{|f(y)|\,|g(x)|}{(x+y)\cos\theta}\,dx\,dy
 \leq C(\theta)\|f\|_2\,\|g\|_2.
\ea
$$
In the last line we refer to the classical (``Hankel-type Hilbert'') inequality for the quadratic
form with kernel $(x+y)^{-1}$
\cite[\S~316]{HLP}.
\eop

\medskip\noindent
{\bf Corollary}. {\it The operator $\Lr$ is bounded in $L^2(\RR_+)$.}
\@ (Take $f=g$ in (\ref{tripleint}).)

\medskip\noindent
Part (a) of the Theorem is thus proved.

\medskip
The main object of study in the proof of part (b)
will be the integral operator defined by the formula
\beq{defS}
 \Sta f(x)=\int_{0}^\infty \frac{f(y)}{x\bom+y\om}\,dy.
\eeq
For $f\in L^2(\RR_+)$ and any $x>0$ the integral converges absolutely by Cauchy-Schwarz.

\begin{lem}
\label{lls}
The operator $\Sta$ is bounded in $L^2(\RR_+)$, and
the formula (\ref{qf}) holds.
\end{lem}

\pf
Lemma \ref{l2absconv} justifies application of  Fubini's theorem to the triple
integral
in the right-hand of (\ref{ltripli}).
Integration with respect to $\rho$ with $x$ and $y$ fixed yields the kernel of the operator $\Sta$
\beq{ints}
\int_0^\infty
e^{-x\rho\bom-y\rho\om}\,d\rho=\f{1}{x\bom+y\om}.
\eeq
Hence $(\Lr f,\Lr g)=(\Sta f,g)$ for all $f,g\in L^2(\RR_+)$.
Combining this with (\ref{tripleint}),
we immediately obtain both assertions of the Lemma.
\eop

\medskip
Our next goal is to find the norm of the operator $\Sta$
precisely. This will be done in Lemmas \ref{ef}--\ref{maxl}.

\medskip
Introduce a family of generalized eigenfunctions which will be used in the spectral decomposition
of the operator $\Sta$:
\beq{psi}
 \psi_\tau(x)=\frac{1}{\sqrt{2\pi}}\;x^{-\frac12+i\tau},\quad x>0,\quad \tau\in\RR.
\eeq
We call $\psi_\tau$ a generalized eigenfunction, since it formally satisfies the equation
$\Sta\psi_\tau=\lam_\tau\psi_\tau$ (see next Lemma); however, $\psi_\tau\notin L^2(\RR_+)$.

\begin{lem}
\label{ef}
 For any $\tau\in \RR$
 \beq{lams}
  \int_0^\infty \f{\psi_\tau(y)}{x\bom+y\om}\,dy=\lam_\tau(\theta)\,\psi_\tau(x),
 \eeq
where
\beq{evalam}
\lam_\tau(\theta)=\f{\pi e^{2\theta\tau}}{\cosh \pi\tau}.
\eeq
\end{lem}

\pf Set $s=-\f{1}{2}+i\tau$. We have (making the change of variable $y\mapsto yx$):
$$
 \int_0^\infty \frac{y^s\,dy}{x\bom+y\om}=
\int_0^\infty \frac{(xy)^s\,x\,dy}{x\bom+(xy)\om}=x^s\int_0^\infty\frac{y^s\,dy}{\bom+\om y}.
$$
The integral converges absolutely, since $\Re s=-1/2$. We obtained (\ref{lams}) with

$$
\lam_\tau(\theta)=\int_0^\infty\frac{y^s\,dy}{\bom+\om y}.
$$
To prove (\ref{evalam}), consider the integral
$$
 I(p)=\int_0^\infty\frac{y^s\,dy}{y+p}\,dy.
$$
If $p>0$, then the substitution $y\mapsto py$ yields
$$
I(p)=B(s+1,-s)p^s,
$$
where $B(\cdot,\cdot)$ is Euler's
Beta function.
By analytic continuation this formula remains valid for all $p\notin [0,\infty)$.
Using the evaluation
$$
B(s+1,-s)=\f{\pi}{\sin\pi(s+1)}=\f{\pi}{\sin(\pi/2+i\pi\tau)}=\f{\pi}{\cosh\pi\tau},
$$
we finally get
$$
\ba{l}\dst
\lam_\tau(\theta)=\bom I(\bom^2)
=\f{\bom \pi}{\cosh\pi\tau}\,(\bom)^{2(-1/2+i\tau)}
\dst
=\f{\pi}{\cosh\pi\tau} e^{(-2i\theta)(i\tau)},
\ea
$$
which simplifies down to (\ref{evalam}).
\eop

\medskip
We are now ready to describe the spectral decomposition
for the operator $\Sta$ explicitly.

\begin{lem}
\label{lemplancherel}
(a)
For any $f\in L^2(\RR_+)$ with compact support in $(0,\infty)$ define
\beq{gencoef}
 \tlf (\tau)=\int_0^\infty f(x)\,\overline{\psi_\tau(x)}\,dx.
\eeq
Then $\tlf\in L^2(\RR)$ and the Plancherel-like identity holds
\beq{genplancherel}
\int_{-\infty}^\infty |\tlf(\tau)|^2\,d\tau=\int_0^\infty |f(x)|^2\,dx.
\eeq
The operator $U:\, f(x)\mapsto \tlf(\tau)$ extends to an unitary operator from $L^2(\RR_+)$
to $L^2(\RR)$
and the inversion formula $f=U^{-1}\tlf$ in the explicit form is
\beq{genfourier}
f(x)=L^2\mbox{-}\lim_{T\to\infty}
\int\limits_{-T}^T \tlf(\tau)\, \psi_\tau(x)\,d\tau.
\eeq
(b)
The action of the operator $\Sta$ on $f\in L^2(\RR_+)$ can be described as follows:
\beq{spdecomp}
 \Sta f= U^{-1}\circ M_{\theta}\circ U\, f,
\eeq
where $M_{\theta}$ is the multiplication operator in $L^2(\RR)$,
\beq{SM}
 M_{\theta}\tlf(\tau)=\lambda_\tau(\theta)\tlf(\tau).
\eeq
\end{lem}

\pf
(a)
Let
\beq{xtov}
x=e^v,
\qquad
g(v)=f(e^v) e^{v/2}.
\eeq
Then
$$
\int_{-\infty}^\infty|g(v)|^2\,dv=\int_{-\infty}^\infty |f(e^v)|^2\, de^v=\int_0^\infty |f(x)|^2\,dx,
$$
so $g\in L^2(\RR)$ and $\|g\|_{L^2(\RR)}=\|f\|_{L^2(\RR_+)}$. Now,
by (\ref{psi}) and (\ref{gencoef})
$$
\tlf(\tau)=\frac{1}{\sqrt{2\pi}}\,\int_{-\infty}^\infty f(e^v) e^{-v/2} e^{-iv\tau}\, de^{v}
=\frac{1}{\sqrt{2\pi}}\,\int_{-\infty}^\infty g(v) e^{-iv\tau}\, dv.
$$
By Plancherel's identity for the Fourier transform, $\|\tlf\|_{L^2(\RR)}=\|g\|_{L^2(\RR)}$.

Hence the transformations $f(x)\mapsto g(v)\mapsto \tlf(\tau)$ are isometries. Therefore
we can extend $U$ to an isometry $L^2(\RR_+)\to L^2(\RR)$.

Expressing $g(v)$ as the inverse Fourier transform of $\tlf(\tau)$ and then expressing
$f(x)$ via $g(v)$ by (\ref{xtov})
yields (\ref{genfourier}). We conclude that $U$ is surjective, hence it is unitary.

\smallskip
(b)
The formal calculation using (\ref{defS}) and (\ref{lams}) goes:
\beq{changeord}
\ba{rl}
\dst
 (U\Sta f)(\tau)&=\dst
 \int_0^\infty\left(\int_0^\infty
  \f{f(y)}{x\bom+y\om}\,dy\right)\,\overline{\psi_\tau(x)}\,dx
\\[3ex]
&=\dst
\int_0^\infty\int_0^\infty\overline{\left(
  \f{\psi_\tau(x)}{ y\bom+x\om }
\right)}\,dx \,f(y)\,dy
\\[3ex]
&=\dst
\int_0^\infty
  \overline{\lambda_{\tau}(\theta)}\;\overline{\psi_\tau(y)}\,f(y)\,dy
\\[3ex]
&=\dst
\lambda_{\tau}(\theta)
\tlf(\tau)
\;=(M_{\theta} U) f(\tau).
\ea
\eeq
In the last line we used (\ref{SM}) and the fact that $\lambda_{\tau}(\theta)\in\RR$,
see (\ref{evalam}).

We will justify the change of order of integration assuming that $f(x)$
has compact support $[a,b]\subset(0,\infty)$. In this case there exists $C_1>0$
(depending on $a$ and $\theta$)
such that $|x\bom+y\om|\geq C_1(1+x)$ for all $x>0$ and all $y\geq a$.

Recalling the definition (\ref{psi}) of $\psi_\tau(x)$, we
conclude that the double integral corresponding to the iterated integrals in
(\ref{changeord})
converges absolutely:
$$
\ba{ll}& \dst
\int_0^\infty\int_0^\infty
  \left|\f{f(y)\,\overline{\psi_\tau(x)}}{x\bom+y\om}\right|\,dx\,dy
\\[3ex]
\leq & \dst
\left(\int_a^b |f(y)|\,dy\right)\,\left(\int_0^\infty \frac{(2\pi)^{-1/2}\,C_1^{-1}}{(1+x)\sqrt{x}}\,dx
\right)
\\[3ex]
\leq & \dst
(b-a)^{1/2}\,\|f\|_2\, C_2\;<\infty.
\ea
$$
By Fubini's theorem, the calculation (\ref{changeord}) is justified, so the formula
$U\Sta f=M_{\theta}Uf$ holds for $f\in L^2$ with compact support in $(0,\infty)$.
Applying $U^{-1}$ on the left in both sides of this formula yields (\ref{spdecomp}) for such $f$.
All operators here are bounded and the set of  compactly supported
functions is dense in $L^2(\RR_+)$. Therefore (\ref{spdecomp})
is valid for all $f\in L^2(\RR_+)$.
\eop

\medskip\noindent
{\bf Corollary}. By (\ref{SM}), the norm of the operator $\Sta$ in $L^2(\RR_+)$ is
\beq{nS}
\|\Sta\|_{2,2}=\|M_{\theta}\|_{2,2}=\sup_{\tau\in\RR} \lambda_\tau(\theta).
\eeq
(Since $\lambda_\tau(\theta)>0$ for all $\theta$ and $\tau$, the absolute value sign is omitted.)

\medskip\noindent
{\bf Remark}. The unitary operator $U$ constructed in Lemma \ref{lemplancherel} does not depend on $\theta$
(the primary reason being that the generalized eigenfunctions (\ref{psi}) are common for
all operators $\Sta$).
As a consequence,
{\it the operators $\St_{\theta}$ and $\St_{\theta'}$ commute for all $\theta,\theta'\in(-\pi/2,\pi/2)$}.
Note however that operators $\Lr$ with different values of $\theta$ do not commute. For instance, $\Lt_{-\theta}\Lr=\Lr^*\Lr=\Sta$,
while $\Lr\Lt_{-\theta}=\St_{-\theta}\neq\Sta$ if $\theta\neq 0$, because of (\ref{spdecomp}), (\ref{SM}) and
the inequality $\lambda_\tau(\theta)\neq\lambda_\tau(-\theta)$.

\medskip
To complete the proof of the formula (\ref{Ltnorm}) and hence of
part (b) of the Theorem, it remains, in view of (\ref{normLS}) and (\ref{nS}),
to evaluate $\sup_{\tau\in\RR}\lambda_\tau(\theta)$.

\begin{lem}
\label{maxl}
For $\lambda_\tau(\theta)$ as defined in (\ref{evalam})
and $K_1(\theta)$ as defined in (\ref{Ltnorm}) we have
\beq{maxlam}
\max_{\tau\in\RR} \lam_\tau(\theta)= (K_1(\theta))^2,
\qquad|\theta|<\pi/2.
\eeq
\end{lem}

\pf
We may assume that $\theta\geq 0$.
Logarithmic differentiation of $\lam_\tau(\theta)$ with respect to $\tau$
yields the condition of extremum:
$\,2\theta-\pi\tanh\pi\tau=0$.
Using the notation $\mu=2|\theta|/\pi$ as in (\ref{Ltnorm}), rewrite this condition as
$\tanh\pi\tau=\mu$.
Consequently, at the point of extremum (which is easily seen to be maximum)
$$
 e^{2\pi\tau}=\frac{1+\mu}{1-\mu}
$$
and
$$
\cosh\pi\tau=(1-\mu^2)^{-1/2}.
$$
We obtain
$$
\max_{\tau\in\RR}
\lambda_\tau(\theta)=\pi\,\left(\frac{1+\mu}{1-\mu}\right)^{\theta/\pi}(1-\mu^2)^{1/2},
$$
which is equivalent to (\ref{maxlam}).
\eop

\medskip
Part (b) of the Theorem is thus proved.

\medskip
To prove part (c), we need to show that $K_1(\theta)$
from (\ref{Ltnorm}) is an increasing function on $[0,\pi/2)$ and to check the formula
(\ref{L2norm}), which also implies left continuity of $K_1(\theta)$ at $\theta=\pi/2$.

This follows from the next, last lemma, where for simplicity
we consider the function
$$
k(\mu)= (K_1(\theta)/\sqrt{\pi})^4=(1-\mu)^{1-\mu}\,(1+\mu)^{1+\mu}.
$$

\begin{lem}
\label{suplam}
The function $k(\mu)$ is increasing in $[0,1)$
and
$$
\lim_{\mu\to 1^-} k(\mu)=4.
$$
\end{lem}

\medskip
\pf
By logarithmic differentiation, $[\ln(x^x)]'=\ln x+1$, hence
$$
 \f{k'(\mu)}{k(\mu)}=-\ln(1-\mu)+\ln(1+\mu)>0.
$$
Finally, $\lim_{\mu\to 1^-}k(\mu)=2^2\,\lim_{x\to 0^+} x^x=4$.
\eop

\medskip
The proof of Theorem \ref{lray2} is complete.


\subsection{Case $p>2$ }
\label{lpneg}

\subsubsection{Why $\Lr:\,L^p\not\to L^{p'}$}
\label{noLp}

It is well known that the Hausdorff-Young theorem for the Fourier transform
does not have an analog for $p>2$ \cite[\S~4.11]{T}.
Moreover, the Fourier transform of an $L^p$ function
with $p>2$ in the sense of distributions  may not even be a locally integrable function \cite[Th.~7.6.6]{ALPDO}.
One explicit example \cite[\S~15.4]{Ed} of a similar kind in the theory of Fourier series is given by
$$
 f(x)=\sum_{n=1}^{\infty} n^{-1/2} e^{ix 2^n},
$$
which is a distribution $\notin L^1(0,2\pi)$.
%

One might expect a different situation for the Laplace transform
and generally for the operators $\Lr$ with $|\theta|<\pi/2$, due to the
exponentially decreasing kernel $e^{-x\rho e^{i\theta}}$.
Such an expectation turns out to be valid only partially.
Lemmas~\ref{lraywl} and \ref{lraywneg} below tell a rather trivial $L^p$ story. Sect.~\ref{lpLorentz} presents a positive result involving Lorentz spaces.

\begin{lem}
\label{lraywl}
If $f\in L^p(\RR_+)$, $p\geq 2$, $|\theta|<\pi/2$, then $\Lr f\in \left(L^2 
+L^{1,\infty}\right)(\RR_+)$.
Moreover, $\Lr f(\rho)$ is continuous at any point $\rho>0$.
\end{lem}

\noindent
{\bf Remark}. The last statement of Lemma shows that $\Lr f$  can be viewed 
as a distribution of order 0 in the space of distributions $\mathcal{D}(\Rplus)\;$ 
\cite[Def.~2.1.1]{ALPDO}. In particular, all integrability troubles of $\Lr f(\rho)$ are attributable
to its limit behaviour at zero or infinity.

\medskip\noindent
\pf
We may assume that $f\geq 0$. Set $f_\infty(t)=\min(f(t),1)$ and write $f=f_\infty+f_2$.
Then $f_\infty\in L^{\infty}$ and $f_2^2(t)\leq f_p^p(t)$, so $f_2\in L^2$. We have
$\Lr f_2\in L^2$  (by Theorem \ref{lray2}) and  
\beq{Linftyweak}
 |\Lr f_\infty (\rho)|\leq\,
 \|f_{\infty}\|_{\infty}\,\int_0^\infty e^{-\rho t\,\cos\theta}\,dt\,=\,
 \f{\|f\|_\infty}{\rho\cos\theta},
\eeq
hence $\Lr f_\infty\in L^{1,\infty}(\RR_+)$.

To prove continuity of $g=\Lr f$ it is enough, again, to separately consider two cases: 
$f\in L^\infty$ and $f\in L^2$. For $f\in L^\infty$ we set 
$f_n(t)=\chi_{(0,n)}(t) f(t)$, where $\chi_I$ is the indicator function of the interval $I$.
Now, $g_n=\Lr f_n$ is easily seen to be continuous for any $n>0$, while
$$
 |g(\rho)-g_n(\rho)|\leq \|f\|_\infty \,\int_n^\infty e^{-t\rho\cos\theta}\,dt
\leq C(\theta,\rho) \frac{\|f\|_\infty}{n}. 
$$
Continuity of $g$ follows by the standard $\eps/3$ trick (cf.~\cite[\S~I.5]{RS}):
$$
 |g(\rho)-g(\rho')|\leq |(g-g_n)(\rho)|+|(g-g_n)(\rho')|+|g_n(\rho)-g_n(\rho')|.
$$
For $f\in L^2$, we may assume $f\geq 0$ and set $f_n(t)=\min(f(t),n)$.
Then $f_n\in L^\infty$, so $g_n=\Lr f_n$ are continuous on $(0,\infty)$ by the previous. By Cauchy-Schwarz,
$$
 |g(\rho)-g_n(\rho)|\leq \|f-f_n\|_2 \,\left(\int_0^\infty e^{-2t\rho\cos\theta}\,dt\right)^{1/2}
\leq C(\theta,\rho) \|f-f_n\|_2. 
$$ 
Since $\|f-f_n\|_2\to 0$ as $n\to\infty$, the continuity of $g(\rho)$ again follows
by the $\eps/3$ trick. 
\eop

\bigskip
\noindent{\bf Remark}.
Although $\Lr:\,L^\infty\to L^{1,\infty}$ and $\Lr:\,L^2\to L^2$, a Marcinkiewicz-type interpolation result is not available. A known subtlety of Marcinkiewicz's theorem
is the required inequality between the exponents of the interpolated spaces, which is violated here: $p=2\not\leq 1=q$. (See e.g.\ \cite{BL}, inequality (7) in \S~1.3 and 
Exercise 9 in \S~1.6.) 

The failure of the method is due to the falsity of the conjecture.    
The truth is, $\Lr$ does not map $\Lp$ to $\Lpp$ for $p>2$
(Lemma~\ref{lraywneg}). This fact 
seems to be not as widely publicized,
even for the standard Laplace transform ($\theta=0$), 
as the corresponding result for the Fourier transform.

The assertion of Lemma~\ref{lraywneg} follows from a result of Bloom \cite[Th.~3.2]{Bloom},
where more general weighted estimates for the Laplace transform are considered.
Our proof offers an explicit counterexample.


\begin{lem}
\label{lraywneg}
$\Lr$ does not map $\Lp$ to $\Lpp$ if $p>2$.
\end{lem}

\pf For simplicity we assume $\theta=0$. For $|\theta|>0$ the same construction works.
Consider the family of functions
$$
 u_\al(t)=t^{\al-1}\,e^{-t}.
$$
Then
$$
 v_\al(x)=\Lt u_\al(x)=\int_0^\infty t^{\al-1}\,e^{-t(1+x)}\,dt=(x+1)^{-\al}\,\Ga(\al).
$$
It is clear that $u_\al\in \Lp(\RR_+)$ if $(\al-1)p>-1$ and $v_\al\in \Lpp(\RR_+)$
if $p'\al>1$. Both conditions are equivalent to $\al>1/p'$, which is assumed in the sequel.

We have
$$
 \|u_\al\|_p^p=\int_0^\infty t^{p(\al-1)}\,e^{-tp}\,dt=p^{-1-p(\al-1)}\,\Ga(p(\al-1)+1),
$$
so
\beq{pnormu}
\|u_\al\|_p=p^{1/p'-\al}\,(\Ga(p(\al-1)+1))^{1/p}.
\eeq
Also,
$$
 \|(x+1)^{-\al}\|_{p'}^{p'}=\int_0^\infty (x+1)^{-p'\al}\,dx=(p'\al-1)^{-1},
$$
so
\beq{ppnormv}
\|v_\al\|_{p'}=\Ga(\al)\,(p'\al-1)^{-1/p'}.
\eeq

Set $\al=1/p'+\eps$. Then $p(\al-1)+1=p(\al-1+1/p)=p\eps$.
Now let $\eps\to 0^+$. The asymptotics of the norm (\ref{pnormu})
is
$$
 \|u_{1/p'+\eps}\|_{p}\sim p^0\,(\Ga(p\eps))^{1/p}\sim C_1\eps^{-1/p}, \quad C_1=C_1(p)>0.
$$
At the same time, the asymptotics of the norm (\ref{ppnormv}) is
$$
 \|v_{1/p'+\eps}\|_{p'}\sim \Ga(1/p')\,(p'\eps)^{-1/p'}\sim C_2(p)\eps^{-1/p'}, \quad C_2=C_2(p)>0.
$$
Therefore
$$
 \frac{\|v_{1/p'+\eps}\|_{p'}}{\|u_{1/p'+\eps}\|_{p}}
 \sim C_3(p)\,\eps^{1/p-1/p'}\to \infty,
$$
since $1/p-1/p'<0$.

\smallskip
To finish the proof, a standard functional-analytic argument is used.
Consider the sequence
$$
\ph_n=\f{u_{1/p'+1/n}}{(\|u_{1/p'+1/n}\|_{p}\,\|v_{1/p'+1/n}\|_{p'})^{1/2}}\;\in L^p.
$$
and let $\psi_n=\Lt \ph_n$.
Then
$\|\ph_n\|_p\to 0$ and $\|\psi_n\|_{p'}\to\infty$ as $n\to\infty$.

Suppose, contrary to the assertion of the Lemma, that $\Lt w\in L^{p'}$ for any $w\in L^p$.
Then $\psi_n\to 0$ weakly in $L^{p'}$, since for any $w\in L^p$
we have
$$
\ba{ll}\dst
 |\langle \psi_n, w\rangle|&\equiv
 \dst \left|\int_0^\infty (\Lt \ph_n)(x)\; w(x)\,dx\right|
\\[3ex] &\dst
 =\left|\int_0^\infty \ph_n(t)\;(\Lt w)(t)\,dt\right|\leq \|\ph_n\|_{p}\,\|\Lt w\|_{p'}
 \to 0.
 \ea
$$
As follows from the Banach-Steinhaus uniform boundedness principle,
the sequence $\{\psi_n\}$ is norm bounded in $L^{p'}$, which is not true by construction.

The Lemma is thus proved.
\eop


\subsubsection{Operators $\Lr$ in Lorentz spaces}
\label{lpLorentz}

Let us briefly recall the necessary definitions and facts,
mostly to fix notation.
For details see e.g.\ \cite{BL}, \S~1.3, or \cite{Garling}, especially Th.~10.4.2. 

The {\it decreasing rearrangement}\ of a nonnegative measurable function $f(t)$ defined on a measurable set
$I$ is the function $f^*:\,[0,\infty)\to[0,\infty]$ defined as
$$
 f^*(t)=\inf\{\lambda\geq 0:\; \Df_f(\lambda)\leq t\},
$$
with $\Df_f(\lambda)$ as in (\ref{Df}). (If $\Df_f(\lambda)$ is strictly monotone and continuous, then
$\lambda=f^*(t)$ is just the inverse function to $t=\Df(\lambda)$.)

The Lorentz space $L^{p,\infty}$ coincides with weak $L^p$ space, that is, the space of measurable functions satisfying
the inequality
$$
\|f\|_{L^{p,\infty}(I)}=\sup t^{1/p}\,|f|^*(t)<\infty.
$$
For  $p\geq 1$, $r\in(0,\infty)$, 
the Lorentz space $L^{p,r}(I)$, is the space of measurable functions on $I$ such that
\beq{Lorentzp}
 \|f\|_{L^{p,r}(I)}^r=\int_0^\infty (t^{1/p}\,|f|^*(t))^r\,\frac{dt}{t}\,<\infty.
\eeq
In particular, $L^{p,p}=L^p$ and it is known that 
\beq{Lorentzinclusion}
L^{p,r}\supset L^{p,s}
\quad\mbox{\rm if}\quad r\geq s.
\eeq
The functionals $f\mapsto \|f\|_{L^{p,r}}$ are not norms
(if $r>p$) but quasinorms: they satisfy the quasi-triangle
inequality $\|f+g\|_{L^{p,r}}\leq C(\|f\|_{L^{p,r}}+\|g\|_{L^{p,r}})$
with some $C>1$.

The following is a $p>2$ version of the Hausdorff-Young inequality for the Laplace transform.

\begin{thm}\label{lrayw}
For any  $p\in (2,\infty)$, any $\theta\in (-\pi/2,\pi/2)$ and any $r\in [1,\infty]$
the operator $\Lr$ carries $L^{p,r}(\RR_+)$ into $L^{p',r}(\RR_+)$. 
\end{thm}

\pf
If $f\in L^\infty(\RR_+)$, then
by 
(\ref{Linftyweak})
$\Lr f \in L^{1,\infty}$ and
$$
 \|\Lr\|_{\infty,1*}\leq (\cos\theta)^{-1}.
$$
On the other hand, Theorem \ref{lray2} and Chebyshev's inequality $\|f\|_{p,\infty}\leq\|f\|_p$
for $p=2$ give
$$
 \|\Lr\|_{2,2*}\leq \|\Lr\|_{2,2*}\leq \sqrt{2\pi}.
$$
It remains to apply the Calder\'on-Marcinkiewicz theorem \cite[Th.~5.3.2]{BL}.
\eop

\medskip
\noindent
{\bf Remark 1}.
Two particular cases 
deserve special attention in view of Lemma~\ref{lraywneg}.
The inequality $1<p'<2<p<\infty$ and the inclusion 
(\ref{Lorentzinclusion}) are used here. 
\\[1ex]
(i) For $r=p'$ we obtain: $\Lr:\;L^{pp'}\to L^{p'}$.
The domain is $L_{pp'}\subsetneq L^{p}$. 
\\[1ex]
(ii) For $r=p$ we obtain: $\Lr:\;L^{p}\to L^{p'p}\supsetneq L^{p'}$.

\medskip
\noindent
{\bf Remark 2}.
The theorem does not state that $\Lr$ is surjective in the concerned pairs of spaces.

\medskip
\noindent
{\bf Remark 3}.
We did not attempt to extract quantitave bounds
for the constants in the inequalities $\|\Lr f\|_{L^{p',r}}\leq C(p,r,\theta)\|f\|_{L^{p,r}}$ from
the quantitative forms of the Calder\'on-Marcinkiewicz existing in the literature
(see e.g.\ \cite{Oklander}, \cite{Holmstedt}).
It would be interesting to determine the growth of $C(p,r,\theta)$ as $|\theta|\to \pi/2-0$.



\section{
Well-projected measures
and weighthed estimates for the Poisson, Cauchy, and 
La\-place integrals}

\label{sec:curve}
\newsectnum

\subsection{Well-projected measures and curves}
\label{sec:wpmes}


For a subset $A\subset\CC$ we denote by $A_x$ and $A_y$ its projections on the coordinate axes. If $A$ is Borel, then  $A_x$ and $A_y$ are Lebesgue measurable and we denote by $|A_x|$, $|A_y|$ their respective Lebesgue measures. 

{\defn 
\label{def:wp}
Let $\mu$ be a Borel measure on a domain $D\subset \CC$. We say that $\mu$ is {\em well-projected}\ with projection constants $k_x\geq 0$ and $k_y\geq 0$
if for any Borel set $A\subset D$ the inequlaity
\beq{wpmes}
 \mu(A)\leq k_x |A_x|+k_y |A_y|
\eeq
holds (provided the right-hand side is finite).
}

The class of all such measures will be denoted $\Wp(D,k_x,k_y)$:
$$ 
 \Wp(D,k_x,k_y)=\{\mbox{\rm Borel measures on $D$ satisfying inequality (\ref{wpmes})}  
\}.
$$ 

In the case $D=\CC$ the notation is abbreviated:
$$
 \Wp(k_x,k_y):=\Wp(\CC,k_x,k_y).
$$ 

The classes $\Wp(k_x,k_y)$ are obviously translation invariant but not rotation invariant. In fact, the well-projectedness property is not rotation invariant regardless of the projection constants. 
See Examples~\ref{ex:horcomb},~\ref{ex:cantorsquare} below.

In order to provide a convenient interface between Theorem~\ref{pcurve} (p.~\pageref{pcurve}) and results of Sect.~\ref{sec:ray}, let us accommodate well-projectedness to rotated axes.

{\defn
\label{def:alpha-wp}
Suppose $\xi, \eta$ are new rectangular coordinates related to $x,y$ by
\beq{xieta}
 \xi=x\cos\al+y\sin\al,\qquad \eta=-x\sin\al+y\cos\al.
\eeq
A Borel measure $\mu$ on $\CC$ is {\em $\al$-well-projected}
if the conditions of Definition~\ref{def:wp} are met with
$\xi$ and $\eta$ substituted for $x$ and $y$ respectively
--- 
that is, if
there exist constants $k_\xi, k_\eta\geq 0$ such that
 for any Borel set $A\subset\CC$
\beq{alpha-wpmes}
 \mu(A)\leq k_\xi |A_\xi|+k_\eta |A_\eta|,
\eeq   
where $A_\xi$ and $A_\eta$ are the projections of the set $A$
on the new coordinate axes.
}

\medskip
The class of such measures will be denoted
$\Wp_\al(k_\xi,k_\eta)$. 

\bigskip
Definition~\ref{def:wp}  
merely axiomatizes the property on which our proof of the key Theorem~\ref{pcurve} hinges. In actual applications of that theorem we will be dealing with rectifiable curves carrying arclength measure. It is convenient to narrow down the definition of  well-projectedness accordingly. 

Fix a domain $D\subset \CC$. Let $\ga$ be a rectifiable curve
in $D$ parametrized by arclength $s\in I$, where $I\subset \RR$ is some interval.
Naturally associated with curve $\ga$ is the {\em arclength measure}\ $\mu_\ga=|ds|$ supported on the closure of the set $\ga(I)$.  

We accept common abuse of terminology, by which a parametric curve $\ga$ means both the map $I\stackrel{\ga}{\to} D$ and its range $\ga(I)\subset D$. Where necessary, a clarification will be provided. 

One situation where such a necessity arises occurs when
a curve has returns and multiple passages (going back and forth through the same place), so that it is impossible to retrieve the measure $\mu_\ga$ uniquely given just the {\em set}\ $\ga(I)$ (the geometric curve), cf.\ Remark~4 below
and Example~\ref{ex:foldcurve}.

Curves are relevant for our purposes only as carriers of measure. There is no point to insist on connectedness or on global continuous parametrization.  Hence the following
generalization.

{\defn 
\label{def:compoundcurve}
Let $\mathcal{I}$ be a nonempty finite or countably infinite index set. A family $\ga=\{\ga_n\}_{n\in\mathcal{I}}$ of rectifiable curves is called a {\em compound curve}.
If the ranges of all $\ga_n$ are contained in a domain $D$,
we say that $\ga$ is {\em a compound curve in $D$}.

A compound curve in $D$ can be regarded as a map whose domain is the disjoint union of the domains $I_n$ of the maps $\ga_n$,
$$
\ga:\; \bigsqcup\limits_{n\in\mathcal{I}} I_n\;\to\;D.
$$ 
The associated {\em total arclength measure}\ is defined as
\beq{compoundmes}
 \mu_\ga=\sum_{n\in\mathcal{I}} \mu_{\ga_n}.
\eeq
}

\smallskip
\noindent
{\bf Remarks}.
1. In Definition~\ref{def:compoundcurve} we do not impose any non-intersection conditions on domains or ranges of the different $\ga_n$. As an extreme example, let all $\ga_n$,
$n\in\NN$, be the same identity map of a certain interval $I$.
Then $\ga:\NN\times I\to I$. If $A\subset I$ is a (Lebesgue-)measurable set, then $\mu_\ga(A)$ can assume only two values: $0$ (if $|A|=0$) and $\infty$ (if $|A|>0$). 

\smallskip
2. Talking about compound curves, we do not use the epithet `rectifiable'; rectifiability of the constituents $\ga_n$
is required by definition. Since expressions like ``a rectifiable curve or a compound curve'' can sound as if `compound' might somehow go without rectifiability, we will write ``an {\em ordinary} rectifiable curve or a compound curve'' in such cases.

\smallskip
3. An ordinary rectifiable curve $\ga$ can be identified with the compound curve $\{\ga\}$, the index set being a singleton. 

\smallskip
4. In applications, a compound curve is often seen informally as the {\em union}\ of the (ranges of) its members $\ga_n$. 
If the curve is non-self-overlapping, i.e. if the map $\ga$ is almost everywhere one-to-one, then the set `range of $\ga$' uniquely determines the measure $\mu_\ga$ (naively --- as length, rigorously --- as the one-dimensional Hausdorff measure, see \cite[\S~3.3.4(A)]{EvGa}). 


\smallskip
5. In this paper we do not exploit a weighted version of the total arclength measure (\ref{compoundmes}), although many
meaningful examples of well-projected measures can be produced that way.  

{\defn
Let $\ga$ be an ordinary rectifiable curve or a compound curve. We say
that $\ga$ is well-projected with projection constants $k_x$, $k_y$ if the associated measure $\mu_\ga$ has this property.  
} 

\medskip
We will use slightly abusive shorthands $\ga\in \Wp(D; k_x,k_y)$ or $\ga\in \Wp(k_x,k_y)$ meaning, of course, that $\mu_\ga$ belongs to the given class.

Similarly we will be talking about {\em $\al$-well-projected
curves} (ordinary rectifiable or compound) and write
$\ga\in\Wp_\al(k_\xi,k_\eta)$ instead of a precise but more cumbersome
$\mu_\ga\in\Wp_\al(k_\xi,k_\eta)$. 

\bigskip
\noindent
{\bf Remark}.
The well-projectedness property makes sense not only in the context of one-dimensional rectifiable sets, see Examples \ref{ex:dxdy} and \ref{ex:cantorsquare}.

\subsection{Examples}
\label{sec:wpexamples}

{\bf Catalog of examples}
\begin{enumerate}
\item {Monotone curve}
\item {Graph of a Lipschitz function}
\item {A self-overlapping curve (a), also treated as a compound curve (b)}
\item {Union of monotone curves with projection overlap of finite multiplicity}
\item {Convex curve}
\item {Simple discontinous radial curve}
\item {Horizontal comb --- a non-well-projected curve that becomes well-projected in rotated coordinates}
\item {Boxed curve}
\item {Characteristic function of 
a bounded set in $\CC$}
\item {Cantor square} 
\end{enumerate}

Graph of dependencies:
%
$$
\ba{ccccc}
1 &\longrightarrow& 4 &\longrightarrow& 5 \\
&          &\downarrow && \downarrow \\
6  && 7 &\longrightarrow& 8
\ea
\mbox{
\begin{picture}(1,30)
\put(-45,9){\rotatebox{-45}{$\longrightarrow$}}
\put(-90,-7){\rotatebox{45}{$\longleftarrow$}}
\end{picture}
}
\qquad\qquad
\ba{ccc}
 2 && 3 \\
  && \\ 
 9 &\longrightarrow & 10
\ea
$$

Examples~\ref{ex:moncurve}--\ref{ex:lipcurve}, \ref{ex:convcurve}--
\ref{ex:boxedcurve}
will have continuation in theorems of Sect.~\ref{sec:app}.

\ex{ex:moncurve}  
A rectifiable curve $\ga\subset D$ parametrized by arclength $s$ is {\em monotone}\ if the coordinate
functions $x(s)$, $y(s)$ are monotone. 
We claim that every monotone $\ga$ is in $\Wp(1,1)=\Wp(\CC;1,1)$.
It suffices to prove this for measurable sets $A\subset \ga$ and there we have 
\beq{monwpproof}
 \mu_\ga(A)=\intl_A ds =\intl_A \sqrt{dx^2+dy^2}\,\leq\,\intl_{A}|dx|+\intl_{A}|dy|
\;\;\stackrel{\raisebox{6pt}{$(*)$}}{=}\;\; 
 |A_x|+|A_y|,
\eeq      
as claimed. The step $(*)$ is where monotonicity is relevant.

\ex{ex:lipcurve}  
Let $\ga$ be the graph of a function $y=f(x)$, $x\in I\subset\RR$. Suppose that $f$ is Lipschitz; presicely --- that
\beq{Lipk}
|f(x_1)-f(x_2)|\leq \lam|x_1-x_2|,
\qquad 
\forall x_1,x_2 \in I.
\eeq
Then $\ga\in\Wp(I\times \RR; \sqrt{1+\lam^2},\,0)$.

Indeed, by Rademacher's theorem \cite[p.~81]{EvGa},
$f$ is a.e.\ differentiable. The definition of derivative implies that at points of differentiability
$$
 |dy/dx|\leq \lam.
$$  
Thus, for a measurable subset $A\subset \ga$
$$
 \int_A |ds|\leq \int_{A_x} \sqrt{1+\lam^2}\,dx.
$$ 
Of course, the roles of $x$ and $y$ in this example can be interchanged.

\ex{ex:foldcurve}  
(a) Consider the parametric curve 
$$
\tilde\ga:\quad
x(t)=t-2\sin t, \quad y(t)=0.
$$
The curve $\ga$ is obtained after re-parametrization of
$\tilde\ga$ by arclength. 
Its range $\ga(\RR)
$ is the $x$-axis, $\RR\times 0\subset\CC$, part of which is covered twice.
Let us compute the measure $\mu_\ga$ or rather its nontrivial
part, restriction on $\RR\times 0$.
We have
$$
\ba{rclcl}
\dst
 \min_{t\in(-\pi,\pi]} x(t)
 & = & \dst
 x|_{t=\pi/3} 
 & =& \dst
\f{\pi}{3}-\sqrt{3},
\\[3ex]
\dst
 \max_{t\in(-\pi,\pi]} x(t)
 & = & \dst
 x|_{t=-\pi/3} 
 & =& \dst
-\f{\pi}{3}+\sqrt{3}.
\ea
$$
It is now clear that
$$
\mu_\ga=(1+\chi_J(x)) \,dx,
$$
where $\chi_J(x)$ is the characteristic function of the set
$$
\ba{c}
\dst
 J=\bigcup\nolimits_{n\in\ZZ} J_n,
\\[2ex]
 \dst
 J_n=\,\left[\f{\pi}{3}-\sqrt{3}+2\pi n,\;
-\f{\pi}{3}+\sqrt{3}+2\pi n\right].
\ea
$$
Since $1+\chi_J(x)\leq 2$, we have $\mu_\ga\in \Wp(2,0)$.

\smallskip
(b) The same measure can be obtained through a different interpretation of $\ga$. Namely, consider $\ga$ as a compound curve,
$
 \ga=\{\ga_n\}_{n\in\ZZ},
$
where (ignoring the trivial $y$-component)
$$
 \ga_n=
\begin{cases}
J_{n-1},& n>0 \\[1ex]
\RR, & n=0\\[1ex]
J_n,& n<0.
\end{cases}
$$
and each $\ga_n$ is parametrized by the identity map.

\ex{ex:monunion}  
This is a straightforward but useful generalization of Example~\ref{ex:moncurve} in the spirit of part (b) of 
Example~\ref{ex:foldcurve}.

Consider a compound curve $\ga=\{\ga_n\}$,
where the member curves $\ga_n\subset D$ are monotone and have the following {\em finite projection multiplicity}\ property. We require that for any $x\in \RR$ 
there be at most $P_1$ member curves whose $x$-projections contain $x$:
$$ 
 \#\{n:\; x\in (\ga_n)_x\}\leq P_1,\qquad \forall x\in \RR. 
$$  
A similar requirement is imposed on $y$-projections:
$$ 
 \#\{n:\; y\in (\ga_n)_y\}\leq P_2,\qquad \forall y\in \RR. 
$$  
Then $\ga\in \Wp(D; P_1, P_2)$.
Indeed, the step $(*)$ of the chain \num{monwpproof} is now changed to 
$$
 \intl_A |dx|+\intl_A |dy| \,\leq\, P_1 |A_x|+P_2 |A_y|. 
$$

\ex{ex:convcurve}  A closed curve $\ga\subset\CC$ is {\em convex}\ if it is the boundary of a convex set. 
Any closed convex curve is the union of at most 4 monotone curves. Hence by Example~\ref{ex:monunion} the arclength measure of any closed convex curve is well-projected with $k_x=k_y=4$.
 
\begin{figure}
\centerline{
\begin{picture}(175,300)
\put(0,0){\includegraphics{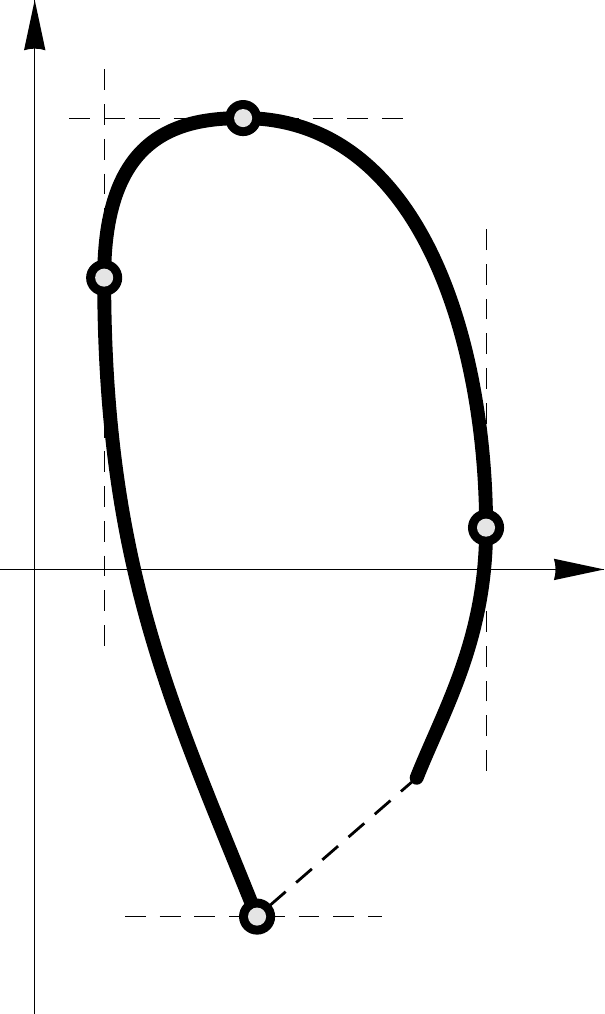}}
\put(164,138){$x$}
\put(20,283){$y$}
\end{picture}
}
\caption{A convex curve 
and its four monotone parts (Example~\ref{ex:convcurve})}
\label{fig:convcurve}
\end{figure}

 A curve $\ga: [a,b]\to\RR$ that is not closed is called convex if the closed curve obtained by joining  the points $\ga(a)$ and $\ga(b)$ with straight line is convex (Fig.~\ref{fig:convcurve}). 

If the domain $I$ of the map $\ga$ is not compact, 
let us exhaust $I$ by compact intervals: 
$I\supset\cdots\supset I_3\supset I_2\supset I_1$.
If each restricted curve $\ga|_{I_n}$ is convex, we say that 
$\ga$ is convex.

It is easy to see that the procedures of closure and compact exhaustion do not affect the projection multiplicity bound 4.
Thus the inclusion
$\ga\in \Wp(D,4,4)$ holds true for any convex curve in $D$,
closed or not.

\ex{ex:maxang}
This example will be referred to in the formulation and proof of 
Theorem~\ref{restrangmax}.

Let $\ph\in(0,\pi/2)$. Consider the sector
$\Phi_\ph=\{z \mid 0<\arg z<\ph\}$. 
Suppose $\al$ is a {\em simple function} defined
on an interval $I\subset \Rplus$ with values in $(0,\ph)$, that is, $\al$ is piecewise constant and assumes only a finite number of values. More specifically,
let $0\leq r_0<r_1<\dots<r_{N-1}<r_N\leq\infty$, $\;I=[r_0,r_N)$,
and $\al(r)=\al_n$ for $r_{n-1}\leq r<r_{n}$, $n=1,2,\dots,N$.
  
Define the corresponding {\em simple discontinuous radial curve}: it is the compound curve
 $\ga=\{\ga_n\}_{n=1}^N$, where $\ga_n$ is the segment
(a ray if $n=N$ and $r_N=\infty$)
connecting the points 
$r_{n-1} e^{i\al_n}$ and $r_{n} e^{i\al_n}$.

Every simple discontinuous radial curve is, of course, well-projected (being a subset of the union of $N$ rays). However there is no uniform estimate; no class $\Wp(\Phi_\ph,k_x,k_y)$ contains all such curves.

Indeed, let us fix $x>0$ and $\psi\in(0,\ph)$.
Take a big integer $M$ and pick $M$ distinct points on the
vertical interval $\Re z=x$, $0<\Im z<x\tan\psi$.
From each of these point issue a tiny radial segment of
length $\eps$. 
If $\eps$ is small enough, there is no radial overlap between these segments and their union $A$ is a subset of some simple discontinuous radial curve. We have $|A|=M\eps$.
At the same time,
$|A_x|<\eps=|A|/M$ and $|A_y|<|A|\sin\psi$. Since $M$ 
can be arbitrarily large and $\psi$
arbitrarily small, an uniform estimate of the form
$|A|<k_x|A_x|+k_y|A_y|$ does not exist. 

\smallskip
We will however prove a {\em restricted uniform estimate}.

Given a sector $\Phi_\ph$ as above, a positive integer $\nu$ and $c>1$, introduce the {\em class}\ $\mathfrak{R}(\ph,\nu,c)$ of simple discontinuous radial curves in $\Phi$ that possess the property
\beq{rgeom}
r_{n+\nu}\geq c^\nu r_{n}\qquad
\mbox{\rm provided}\quad
n+\nu\leq N.
\eeq


Intuitively, this assumption requires the curve to have long radial runs (of length comparable to the current value of $r$) between angular jumps.

\begin{lem}
\label{lem:wprad}
All simple discontinuous radial curves in the sector $\Phi_\ph$ that satisfy condition (\ref{rgeom}) are well-projected with uniformly bounded projection constants: we have the inclusion  
\beq{wpradclass}
\mathfrak{R}(\ph,\nu,c)\subset \Wp(\Phi_\ph;\kappa,0),
\eeq
where
\beq{wpradk}
\kappa=\left(\f{\log\sec\ph}{\log c} +\nu+2\right)\,\sec\ph.
\eeq
\end{lem}

 
\pf
Let $\ga\in \mathfrak{R}(\ph,\nu,c)$ consist of
the segments $\ga_n=[r_{n-1}e^{i\al_n},\,r_n e^{i\al_n})$.
Consider their $x$-projections
$J_n=(\ga_n)_x=[r_{n-1}\cos\al_n,\,r_n\cos\al_n)$.
We want to obtain an upper bound for the multiplicity of the $x$-projections
$$
 P_1=\sup_{x>0}\,\#\{n\mid J_n\ni x\}.
$$ 
(cf.~notation in Example~\ref{ex:monunion}).
Fix $x>0$. Let 
$$
\ba{l}
  m=\min\{n\mid J_n\ni x\}, \qquad 
  M=\max\{n\mid J_n\ni x\}.
\ea
$$
Ignoring the trivial case $P_1=1$, we may assume that $M>m$. 
Then 
$$
 r_m> r_m\cos\al_m >\;\; x\;\; 
\geq r_{M-1}\cos\al_M
> r_{M-1}\cos\ph.
$$
Let $h$ be the positive integer such that 
\beq{hMang}
m+\nu h> M-1\geq m+\nu(h-1). 
\eeq
Then by (\ref{rgeom})
$$
 r_{M-1}\geq c^{\nu(h-1)} r_m.
$$
Comparing to the previous inequality, we see that
$$
 \sec\ph > \f{r_{M-1}}{r_m} \geq c^{\nu(h-1)}.
$$
Therefore
$$
 \nu(h-1)< \f{\log\sec\ph}{\log c}.
$$
Using the left inequality (\ref{hMang}), we get
\beq{radestPx}
 M-m+1 \;<\;\nu h+2 = \nu(h-1)+\nu+2\;< 
\f{\log\sec\ph}{\log c} +\nu+2.
\eeq
The right-hand side does not depends on $x$, hence it can serve as an upper bound for $P_1=\sup_x (M-m+1)$.

Unlike in Example~\ref{ex:monunion}, we don't have an estimate for the multiplicity $P_2$ of $y$-projections; instead we do have bounded slope: $ds/dx<\sec\ph$ on every segment $\ga_n$.
Hence
$$
d\mu_\ga\leq P_1\sec\ph\;dx.
$$
Substituting the estimate for $P_1$ from the r.h.s.\ of (\ref{radestPx}), we obtain (\ref{wpradclass})--(\ref{wpradk}).
\eop

\bigskip
In the next examples we consider $\al$-well-projected measures and curves
(recall Definition~\ref{def:alpha-wp}). 

\ex{ex:horcomb}
A {\em horizontal comb}\ is a compound curve of the form
$$
\ga=\{\ga_n\}_{n\in\mathfrak{I}},
\;\quad
\ga_n=[a_n, a_n+ib_n], \quad \mathfrak{I}=\NN\;\;\mbox{\rm or}\;\;\ZZ.
$$
The segments $\ga_n$ will be referred to as ``teeth''. 
We assume that $a_n>a_m$ when $n>m$, and that all $b_n>0$.
The index set $\NN$ corresponds to a ``semi-infinite'' comb, 
and the index set $\ZZ$ yields a comb with infinitely many teeth in both directions.

Similarly, a {\em vertical comb}\ consists of horizontally
oriented teeth 
$$
\ga_n=[ia_n, ia_n+b_n],\quad n\in\NN\;\mbox{\rm or}\;\ZZ.
$$

The combs are never well-projected.%
\footnote
{Take an arbitrarily large $N$ and
  let $\be=\min(b_1,\dots,b_N)$. If $\pi_y$ denotes
a projection on the $y$-axis, then 
$\mu_\ga(\pi_y^{-1}([0,\be])\geq N\be$, making 
(\ref{wpmes}) impossible.
} 
In a {\em rotated}\ coordinate system, however, well-projectedness with {\em uniformly bounded}\ projection constants takes place for some classes of combs.

Given a positive integer $\nu$ and $c>0$, introduce the {\em class}\ $\mathfrak{C}(\nu,c)$ of semi-infinite horizontal combs that possess the properties
\beq{moncomb}
b_1\leq b_2\leq b_3\leq\dots
\eeq
and
\beq{hcombc}
\f{a_{n+\nu}-a_{n}}{k} \geq  c b_n,
\quad
\forall n\geq 1.
\eeq
The latter property means that the teeth must not be very long compared to the spacings between them (averaged over $\nu$ successive intervals). The role of condition (\ref{hcombc}) here is similar to the role of condition (\ref{rgeom}) in Example~\ref{ex:maxang}.

\begin{lem}
\label{lem:wpcomb}
For any $\al\in(0,\pi/2)$, all semi-infinite combs that satisfy condition (\ref{hcombc}) are $\al$-well-projected with uniformly bounded projection constants: we have the inclusions  
\beq{wpcombclass}
\mathfrak{C}(\nu,c)\subset \Wp_\al(\kappa_1,0)
\qquad
\mbox{\rm and}
\qquad
\mathfrak{C}(\nu,c)\subset \Wp_\al(0,\kappa_2),
\eeq
where
\beq{wpcombk}
\ba{l}
\kappa_1=\dst c^{-1}\,\sec\al\, +\, (\nu+1)\,\csc\al,
\\[1ex]
\kappa_2=\dst c^{-1}\,\csc\al\, +\, (\nu+1)\,\sec\al.
\ea
\eeq
\end{lem}

\pf 
Structurally this proof repeats the proof of
Lemma~\ref{lem:wprad}.

\smallskip

Let $\xi$ and $\eta$ be the new coordinates defined as in
(\ref{xieta}). 
Suppose $\ga\in \mathfrak{C}(\nu,c)$ consists of
the segments (teeth) $\ga_n=[a_n,\,a_n+ i b_n]$.
Consider their projections on the $\xi$ axis,
\beq{combJn}
J_n=(\ga_n)_\xi=[a_{n}\cos\al,\,a_n\cos\al+b_n\sin\al].
\eeq
We want to obtain an upper bound for the multiplicity $P_1$ of the $\xi$-projections
$$
 P_1=\sup_{\xi}\,\#\{n\mid J_n\ni \xi\}.
$$

%
Fix $\xi>0$. Let 
\beq{mMcomb}
  m=\min\{n\mid J_n\ni \xi\}, \qquad 
  M=\max\{n\mid J_n\ni \xi\}.
\eeq
Then 
\beq{aMmcomb}
 a_m \geq \xi\sec\al-b_m\tan\al,
\qquad
a_M\leq \xi\sec\al.
\eeq
Therefore, 
\beq{amMdif}
a_{M}-a_m\leq b_m\tan\al. 
\eeq
Compared to the proof of Lemma~\ref{lem:wprad},
we modify the definition of the positive integer $h$ replacing $M-1$ in the inequalities (\ref{hMang}) by $M$:
\beq{hMcomb}
m+\nu h> M\geq m+\nu(h-1). 
\eeq
Then by (\ref{moncomb})--(\ref{hcombc})
\beq{amMdif2}
 a_{M}-a_m\geq c\nu(h-1) b_m.
\eeq
Comparing to (\ref{amMdif}) we see that
\beq{hcombal}
 c\nu(h-1)\leq\tan\al .
\eeq
Using left inequality (\ref{hMcomb}), we get
$$
 M-m+1 \;<\;\nu h+1 = \nu(h-1)+\nu+1\;\leq 
\f{\tan\al}{c} +\nu+1.
$$
The right-hand side does not depends on $\xi$, hence it can serve as an upper bound for $P_1=\sup_{\xi} (M-m+1)$:
\beq{combestPxi}
P_1\leq \f{\tan\al}{c} +\nu+1.
\eeq
The multiplicity $P_2$ of projection of $\ga$ on the $\eta$-axis is estimated similarly:
\beq{combestPeta}
P_2\leq \f{\cot\al}{c} +\nu+1.
\eeq
Since the teeth are vertical, we have $ds=dy$, so 
$$
\left|\f{ds}{d\xi}\right|=\csc\al,
\qquad
\left|\f{ds}{d\eta}\right|=\sec\al.
$$
Combining these estimates with 
(\ref{combestPxi})--(\ref{combestPeta}), we obtain (\ref{wpcombclass})--(\ref{wpcombk}).
\eop

\ex{ex:boxedcurve} 
Let $\{B_n\}$, $n=1,2,\dots$, be a family of non-overlapping rectangles (``boxes'') with sides parallel to the $x$, $y$ axes. Using the notation of
Example~\ref{ex:horcomb}, let the left lower corner and the right upper corner of
$B_n$ be, respectively, $a_n$ and $a_{n+1}+ib_n$. 
 
\begin{figure}
\centerline{
\begin{picture}(250,155)
\put(0,0){\includegraphics{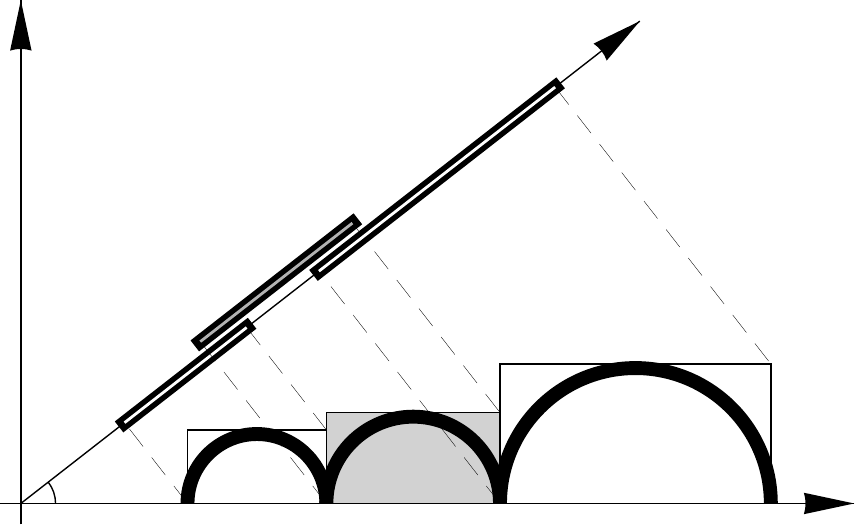}}
\put(234,14){$x$}
\put(13,140){$y$}
\put(163,142){$\xi$}
\put(20,8.5){$\alpha$}
\put(110,12){$B_n$}
\put(64,75){$J_n$}
\end{picture}
}
\caption{Boxed non-Lipschitz curve (Example~\ref{ex:boxedcurve}; Theorem~\ref{thm:boxcurve})
and its $\xi$-projection}
\label{fig:semicircles}
\end{figure}
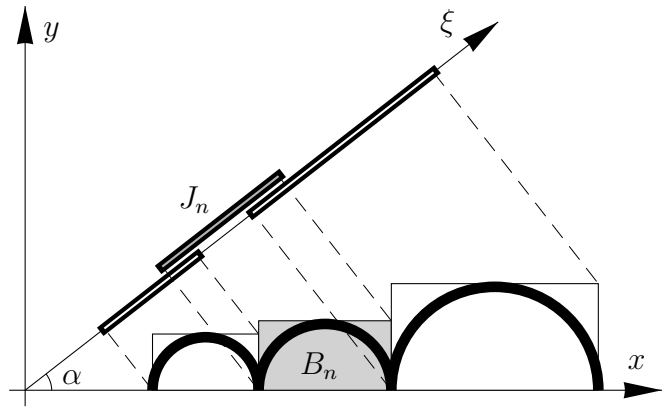

Consider a compound curve $\ga=\{g_n\}$, 
where each $\ga_n$ is {\em a convex curve contained in the
(closed) box $B_n$}\
(Fig.~\ref{fig:semicircles}).
We say that $\ga$ is a {\em boxed curve}\
with parameters $\{a_n\}$, $\{b_n\}$. The horizontal comb of Example~\ref{ex:horcomb}
is a particular case: its teeth are precisely the the left sides of the boxes $B_n$.  

Introduce the {\em class $\mathfrak{B}(\nu,c)$}\ of boxed curves for which the box parameters satisfy the conditions (\ref{moncomb}), (\ref{hcombc}). Thus we require the heights
to be nondecreasing and the aspect ratios `height/widht' bounded on average (over index intervals of length $\nu$).
Similarly to Example~\ref{ex:horcomb},
curves of every such class are well-projected with uniformly
bounded projection constants in a rotated coordinate system. 
 
\begin{lem}
\label{lem:boxproj}
For any $\al\in(0,\pi/2)$ we have the inclusions  
\beq{wpboxclass}
\mathfrak{B}(\nu,c)\subset \Wp_\al(\kappa_1,\kappa_2)
\eeq
where
\beq{wpboxk}
\ba{l}
\kappa_1=\dst 4\left(c^{-1}\,\sec\al\, +\, (\nu+2)\,\csc\al\right),
\\[1ex]
\kappa_2=\dst 4\left(c^{-1}\,\csc\al\, +\, (\nu+2)\,\sec\al\right).
\ea
\eeq
\end{lem} 
 
\pf 
We only mention the changes to be made in the proof
of Lemma~\ref{lem:wpcomb}.

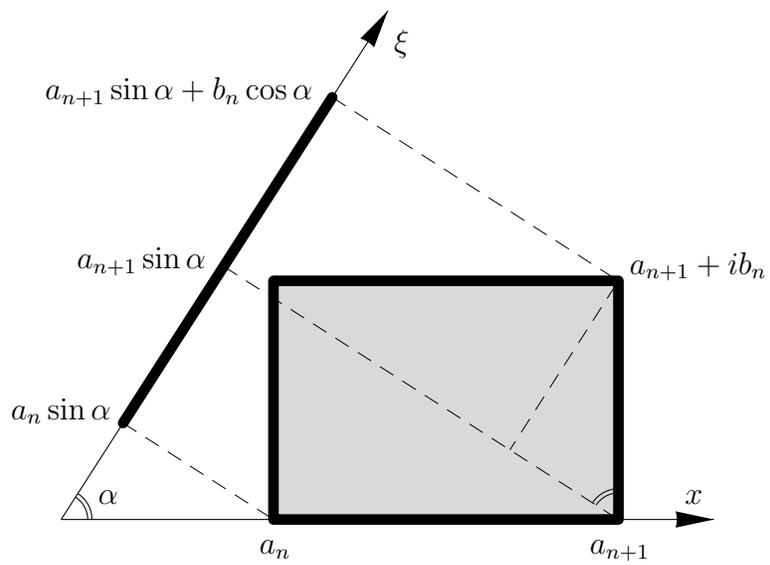
\begin{figure}
\centerline{
\begin{picture}(285,200)
\put(19,10){\includegraphics{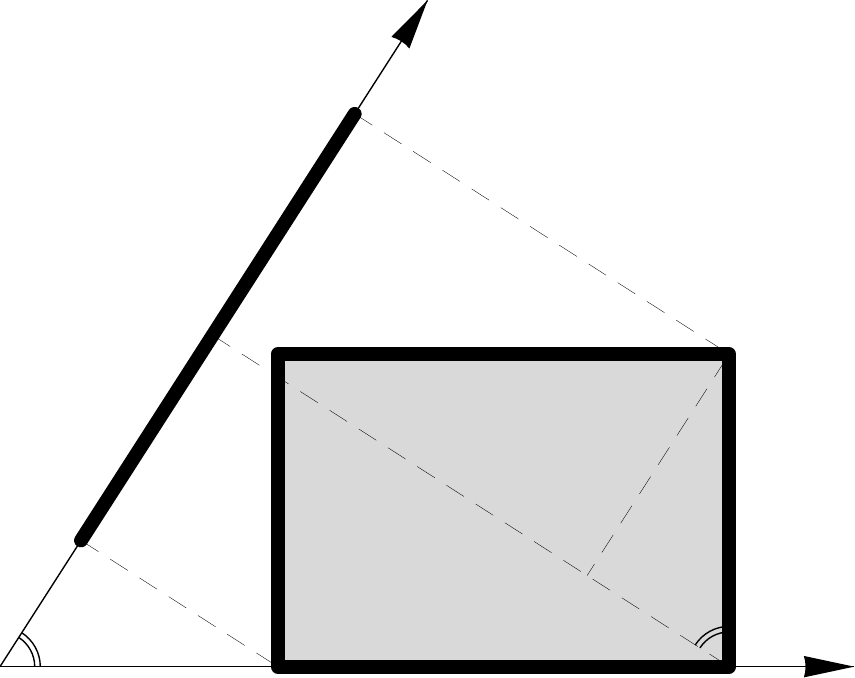}}
\put(33,19){$\al$}
\put(255,19){$x$}
\put(145,190){$\xi$}
\put(94,0){$a_n$}
\put(219,0){$a_{n+1}$}
\put(234,106){$a_{n+1}+ib_n$}
\put(0,51){$a_{n}\sin\al$}
\put(25,108){$a_{n+1}\sin\al$}
\put(13,172){$a_{n+1}\sin\al+b_n\cos\al$}
\end{picture}
}
\caption{Projection of a box onto a ray (Lemma~\ref{lem:boxproj})}
\label{fig:boxproj}
\end{figure}

The segments $J_n$ on the $\xi$-axis (cf.~(\ref{combJn})) are now defined as projections of the boxes $B_n$ (Fig.~\ref{fig:boxproj}),
$$
J_n=(B_n)_\xi=[a_{n}\cos\al,\,a_{n+1}\cos\al+b_n\sin\al].
$$
The inequalities (\ref{aMmcomb}) become
$$
 a_{m+1} \geq \xi\sec\al-b_m\tan\al,
\qquad
a_M\leq \xi\sec\al.
$$
%
Instead of (\ref{amMdif}) we have
\beq{amMdifbox}
a_{M}-a_{m+1}\leq b_m\tan\al.
\eeq
Defining the number $h$, we return from (\ref{hMcomb}) back to (\ref{hMang}).
The inequality (\ref{amMdif2}) becomes
$$
a_M-a_{m+1}\geq c\nu(h-1)b_m.
$$
The inequality (\ref{hcombal}) stays.

The difference between (\ref{hMcomb}) and (\ref{hMang})
results in the extra summand $1$ in the right-hand side of
(\ref{combestPxi}): we now get  
\beq{boxestPxi}
P_1\leq\,\f{\tan\al}{c}+\nu+2.
\eeq
Similarly (by changing $\al\to \pi/2-\al$) we get
\beq{boxestPeta}
P_2\leq\,\f{\cot\al}{c}+\nu+2.
\eeq
Had we assumed the curves $\ga_n$ to be monotone, we would  
have already arrived at the situation of Example~\ref{ex:monunion} relative to the coordinates $\xi$, $\eta$. Since our assumption about $\ga_n$ is convexity,
the extra factor 4 
(cf.~Example~\ref{ex:convcurve}) appears 
in the ultimate projection multiplicity bounds. 

Now we infer (\ref{wpboxclass})--(\ref{wpboxk})
from (\ref{boxestPxi})--(\ref{boxestPeta}) referring to Example~\ref{ex:monunion}.
\eop

\ex{ex:dxdy}  
Let $\chi_G$ be the characteristic function of a 
bounded Borel set $G\subset\CC$.
Define a measure $\mu$ by the formula
$$
 \mu(A)=\int_{A}\chi_G(x,y)\,dx\,dy=\int_{A\cap G}\,dx\,dy,
$$
where $A$ is any Borel set in $\CC$.
In other words, $\mu$ is the restriction of the planar Lebesgue measure $dx\otimes dy$ to the set $G$.
We claim that 
\beq{dxdywp}
\;\mu\in\Wp(t|G_y|,\;(1-t)|G_x|)
\eeq
for any $t\in [0,1]$.

\smallskip
\pf
Without loss of generality we may
assume that $A\subset G$,
so that $A$ is bounded.
Clearly, $A\subset A_x\times A_y$, hence
$$
\mu(A)\leq |A_x|\,|A_y|.
$$
To finish the proof, denote 
$u=|A_x|$, $v=|A_y|$,
$U=|G_x|$, $V=|G_y|$ and use the obvious  inequality
$$
 \f{uv}{UV}\,\leq\; t\,\f{u}{U}+(1-t)\,\f{v}{V}
\qquad \mbox{\rm provided}\quad u\leq U,\;v\leq V.
$$
\eop

\ex{ex:cantorsquare}
Let 
$$
\Can=\bigcap_{n=0}^\infty\Can_n
$$ 
be the standard closed ternary Cantor set. Here $\Can_0=[0,1]$
and the $n$-th iteration set $\Can_n$, $n=1,2,\dots$, is
obtained by removing the middle open interval of length $3^{-n}$ from each of the $2^{n-1}$ intervals consituting $\Can_{n-1}$.

The measure $\mu$ in this example is supported on the Cartesian product $\Can\times\Can$. It will be defined by
a limiting process.

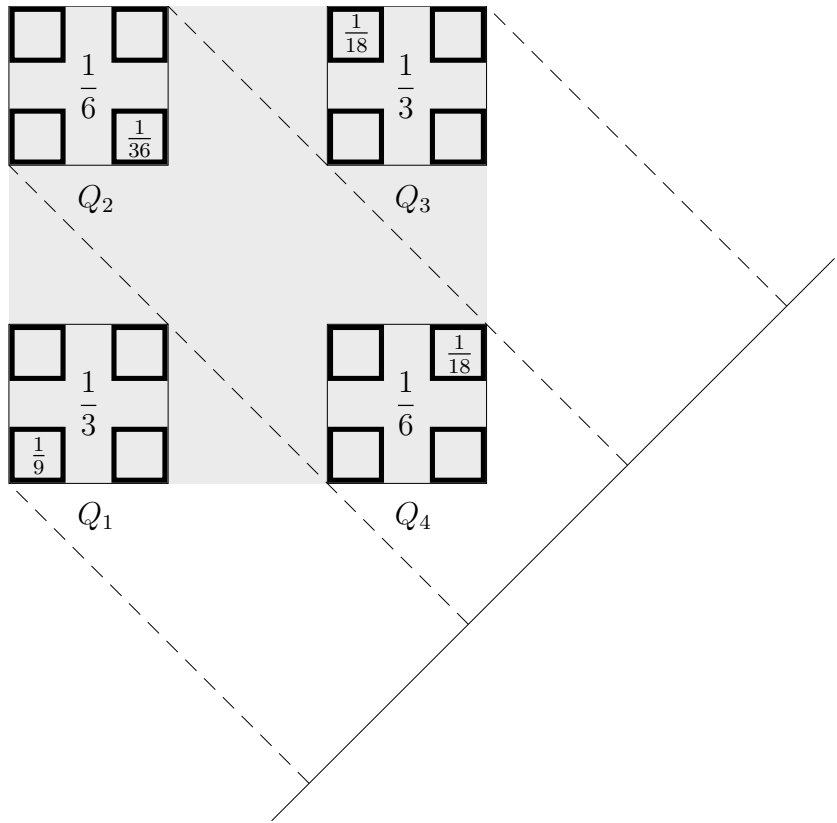
\begin{figure}
\centerline{
\begin{picture}(310,310)
\put(0,0){\includegraphics{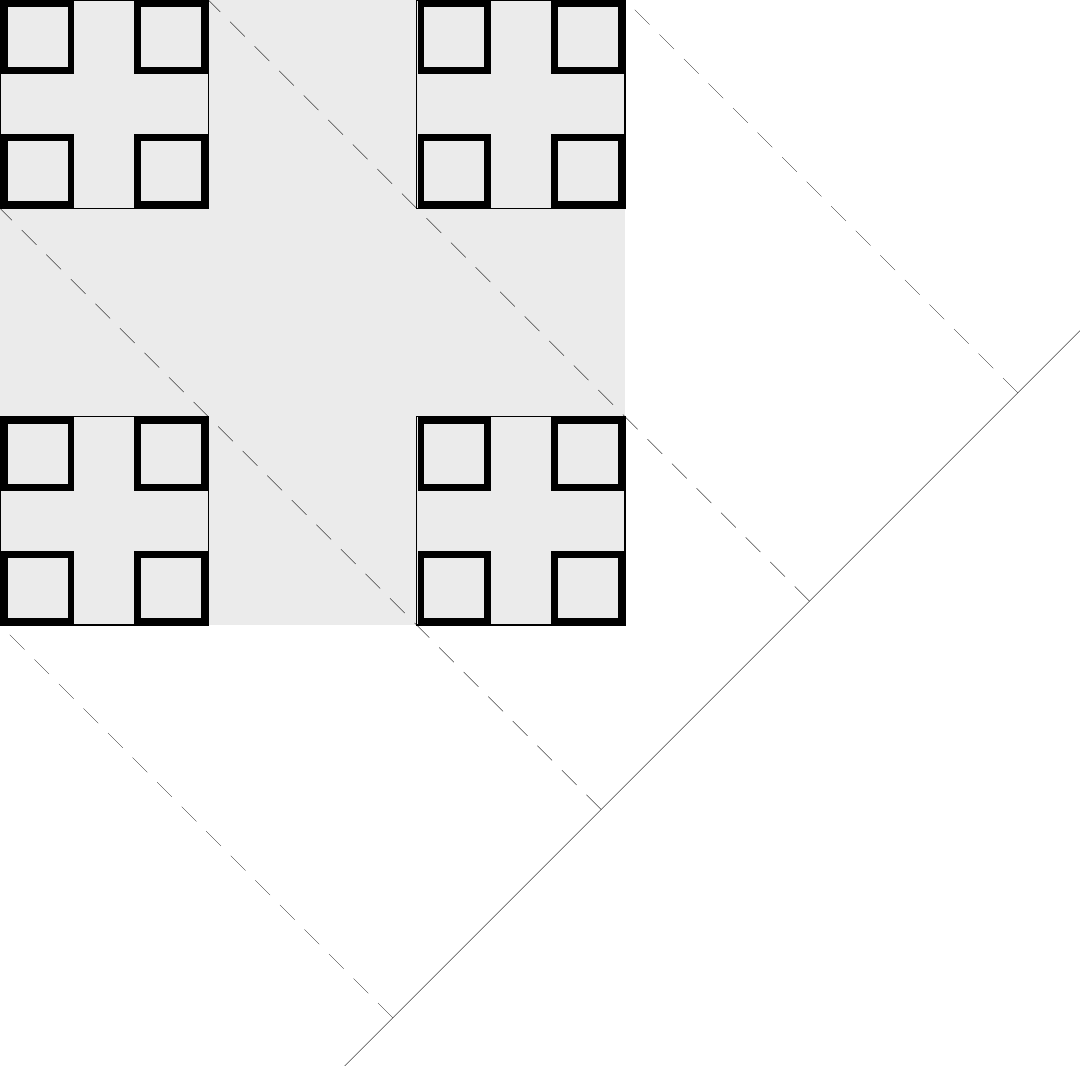}}  
\put(26,154){$\dst\f13$}
\put(146,154){$\dst\f16$}
\put(26,274){$\dst\f16$}
\put(146,274){$\dst\f13$}
\put(8,135){$\textstyle \f{1}{9}$}
\put(165,173.5){$\textstyle \f1{18}$}
\put(44,255){$\textstyle \f1{36}$}
\put(126,295){$\textstyle \f1{18}$}
\put(26,112){$Q_1$}
\put(26,232){$Q_2$}
\put(146,232){$Q_3$}
\put(146,112){$Q_4$}
\end{picture}
}
\caption{Cantor square (Example~\ref{ex:cantorsquare}):
 first and second generation squares}
\label{fig:cantorsquare}
\end{figure}

Let us begin with Lebesgue measure restricted to the unit square, 
$$
\mu_0=dx\otimes dy|_{\Can_0\times\Can_0}.
$$
%
%
Let $Q_1,\dots,Q_4$ be the 
first-generation squares with side length
$1/3$ (Fig.~\ref{fig:cantorsquare}).
Consider the linear maps $T_{j}$
of the unit square $\Can_0$ onto $Q_j$, $j=1,\dots,4$. 
%
%
Let $\mu_0^{(j)}$ be the push-forward of the measure $\mu_0$ 
by the map $T_j$, i.e.\ for any Borel set $A$
$$
 \mu_0^{(j)}(A)= \mu_0(T_j^{-1} (A)).
$$
Since $T_j$ is a composition of the homothety with coefficient $1/3$ and a translation, we have the explicit formula
$$
\mu_0^{(j)}=9\;dx\otimes dy|_{Q_j}.
$$
Define 
$$
 \mu_1=\f13\mu_0^{(1)}+\f16\mu_0^{(2)}+\f13\mu_0^{(3)}+\f16\mu_0^{(4)}.
$$
(The weights are so chosen to make the proof of Lemma~\ref{lem:wpcantor} work.)
Thus $\mu_1(Q_1)=\mu_1(Q_3)=1/3$, 
$\mu_1(Q_2)=\mu_1(Q_4)=1/6$, 
and $\mu_1(\Can_1\times\Can_1)=1$.

Continuing by induction and assuming that the measure $\mu_{n}$
has been constructed, define $\mu_{n}^{(1)},\dots,\mu_{n}^{(4)}$ as the push-forwards 
\beq{mun}
\mu_{n}^{(j)}(A)=\mu_{n}(T_j^{-1}(A)). 
\eeq
Now define
\beq{cantormu}
 \mu_{n+1}=\f13\mu_{n}^{(1)}+\f16\mu_{n}^{(2)}+\f13\mu_{n}^{(3)}+\f16\mu_{n}^{(4)}.
\eeq
Observe that $\mu_n(\Can_1\times\Can_1)=1$ for all $n\geq 0$. More generally, if $Q$ is one of $2^k$ squares of $k$-th generation and $n\geq k$, then $\mu_{n+1}(Q)=\mu_n(Q)$. (Informally speaking, to reveal the difference between $\mu_{n+1}$ and $\mu_n$ one needs a microscope capable to resolve distances as small as $3^{-n-1}$.)

Consequently, if $R\subset\CC$ is a rectangle whose vertex coordinates
are rationals of the form $k/3^m$, $k,m\in\ZZ$,
then the sequence $\{\mu_n(R)\}$ stabilizes: $\mu_{m}(R)=\mu_{m+1}(R)=\dots$, and existence of the limit
$$
 \lim_{n\to\infty}\mu_n(R)=\mu_m(R)
$$
is trivial. We define $\mu(R)$ as this limit.
Using inner and outer approximations, we then extend $\mu$ to all rectangles with sides parallel to the axes, then --- to all open sets, and in the final step of construction --- to the $\sigma$-algebra of all Borel subsets of $\CC$.

\smallskip
Since $\mu(\PCan)=\mu(\CC)=1$ while the $x$- and $y$-projections of $\PCan$ have Lebesgue measure 0, {\em the measure $\mu$ is not well-projected}. We will show, however, that {\em it is well-projected in rotated coordinates}. 
 
\begin{lem}
The measure $\mu$ constructed above belongs to the class
\label{lem:wpcantor}
$
\Wp_{\pi/4}(\sqrt{2},0).
$
\end{lem}

\pf
Let $\xi,\eta$ be the coordinates (\ref{xieta}) with $\al=\pi/4$.
It suffices to prove the estimates
\beq{cantorn}
\mu_n(A)\leq \sqrt{2}\,|A_\xi|
\eeq
for $n=1,2,\dots$, where $A\subset\CC$ is any Borel set. 
 
For $n=0$, the inequality \ref{cantorn} follows from
Example~\ref{ex:dxdy} with $t=1$ (in the rotated axes), since 
$|(\Can_0)_\eta|=\sqrt{2}$.

\unl{Induction step $\,n\,\to\,n+1$}. 
Let $A_j=A\cap Q_j$ ($j=1,\dots,4$).
Combining (\ref{mun}), (\ref{cantorn}), and the fact that
$T_j^{-1}$ expands lengths by a factor of $3$, we get 
$$
\mu_n^{(j)}(A)\leq \sqrt{2}\,\left|\big(T_j^{-1}(A_j)\big)_\xi\right|
\,=\,3\,\sqrt{2}\,\left|(A_j)_\xi\right|.
$$
Therefore by (\ref{cantormu})
\beq{cantormuA}
\mu_{n+1}(A)\leq
\sqrt{2}\left(|(A_1)_\xi|+\f{1}{2}|(A_2)_\xi|
+|(A_3)_\xi|+\f{1}{2}|(A_4)_\xi|\right).
\eeq

The decisive fact: the $\xi$-projections of the sets $A_1$, $A_3$, and $A_2\cup A_4$ do not overlap. Hence
$$
\ba{lcl}
 |A_\xi|&\geq&\dst
 |(A_1)_\xi|+|(A_3)_\xi|
+\max\big(|(A_2)_\xi|,\,|(A_4)_\xi|\big)
\\[2ex]
&\geq&\dst
|(A_1)_\xi|+|(A_3)_\xi|
+\f{1}{2}\left(|(A_2)_\xi|+|(A_4)_\xi|\right).
\ea
$$
Comparing with (\ref{cantormuA}), we conclude that
$$
\mu_{n+1}(A)\leq \sqrt{2}|A_\xi|,
$$
which completes the induction step
and the proof of Lemma.
\eop

%


\subsection{Estimates for the Poisson and Cauchy integrals in $L^p$ spaces with well-projected weights}

\label{sec:cauchy}

Our local goal is to obtain certain weighted $L^p$ estimates of $x$-variable convolutions of $L^p(\RR)$ functions with the Poisson kernel 
$$
 P(x,y)=\frac{1}{\pi}\,\frac{y}{x^2+y^2}
$$
and the Cauchy kernel
$$
 S(x,y)=\frac{i}{\pi}\,\f{1}{x+iy}.
$$
For convenience we will always assume $y>0$. Thus
the complex number 
$$
z=x+iy
$$ 
will belong to the upper half-plane $\HH$. 

Consider first the Poisson integral
\beq{Poissonint}
 \cP u(z)=\int_{\RR} P(x-t,y)\,u(t)\,dt.
\eeq
If $u\in L^\infty(\RR)$, the integral is pointwise defined everywhere in $\HH$ due to the properties of the Poisson kernel:
\beq{propP}
P(x,y)>0,\qquad \int_{\RR} P(x,y)\,dx=1, \quad \forall y>0,
\eeq
and the obvious estimate holds:
\beq{PLinfty}
 |\cP u(z)|\leq\|u\|_{\infty},\qquad \forall z\in\HH.
\qquad
\eeq

Suppose $\mu$ is a well-projected measure in $\HH$.
The weighted spaces $L^p(d\mu)$, $1<p\leq \infty$, and the weighted weak space $L^{1,\infty}(d\mu)$ are defined in an obvious way. For instance, if $v(z)$ is a $\mu$-measurable function in $\HH$, then
$$
 \|v\|_{L^{1,\infty}(d\mu)}=
\sup\limits_{\lam>0}\left(\lam\,\cdot\,\mu\{z\in\HH:\; |v(z)|>\lam\}\right).
$$

\begin{thm}\label{pcurve}
Suppose 
\beq{muwpclass}
\mu\in \Wp(\HH; k_x,k_y).
\eeq
(a) The operator $\cP$ acts continuosly from $L^1(\RR)$ to $L^{1,\infty}(d\mu)$. Quantitatively,
\beq{PL1wmu}
 \|\cP u\|_{L^{1,\infty}(d\mu)}
\leq
 (3 k_x+\pi^{-1} k_y)\,\|u\|_{L^1(\RR)}.
\eeq
(b)
Let $1<p\leq\infty$. 
The operator $\cP$ defined initially on $L^\infty(\RR)\cap L^p(\RR)$ extends to a bounded
operator from $L^p(\RR)$ to $L^p(d\mu)$ with uniform 
norm estimate
\beq{estgamA}
\KPcurve(p,k_x,k_y):=\sup_{\mu\in\Wp(\HH,k_x,k_y)} \|\cP\|_{L^p(\RR)\to L^p(d\mu)}
 \leq 2\left((3 k_x+\pi^{-1} k_y)p'\right)^{1/p}.
\eeq
\end{thm}

\pf
(a) 
Fix $\lam>0$. Consider the set
\newcommand{\Projset}{A} 
\beq{weakAgam}
 \Projset=\Ds_{\cP u}(\lambda)=\{z\in \HH :\, |\cP u(z)|>\lambda\}
\eeq
and estimate the Lebesgue measure of its projections 
$ 
\Projset_x 
$ 
and  
$ 
\Projset_y 
$ 
on the coordinate axes.

\label{pg:projest}
\unl{Estimation of $|\Projset_y|$} is elementary.
Note that $P(x,y)\leq (\pi y)^{-1}$.
Therefore
$|\cP u(x+iy)|\leq \|u\|_1/(\pi y)$.
Consequently, if $y\in \Projset_y$, then $\lambda<\|u\|_{1}/(\pi y)$, hence
\beq{mesSy}
|\Projset_y|\leq \f{\|u\|_{1}}{\pi \lam}.
\eeq
	
\unl{Estimation of $|\Projset_x|$}. Consider the (centered) Hardy-Littlewood maximal function
\cite[\S~2.1]{Grafakos} of $u(t)$,
$$
U(t)=\sup_{a>0}\f{1}{2a}\int_{t-a}^{t+a}|u(t)|\,dt.
$$
Because of the properties \num{propP} we have $|\cP u(x+iy)|\leq U(x)$ (cf.\ \cite[Th.~2.1.10]{Grafakos}
or 
see p.~\pageref{PleqU}).
Therefore $\Projset_x\subset\Ds_U(\lam)$. By the Hardy-Littlewood maximal theorem \cite[Th~2.1.6]{Grafakos} we obtain
\beq{mesSx}
 |\Projset_x| \leq |\Ds_U(\lam)|\leq \f{3\|u\|_{1}}{\lam}.
\eeq
\label{pg:eoppcurve}
Combining the inequalities (\ref{mesSx}), (\ref{mesSy}) and (\ref{wpmes})
we obtain (\ref{PL1wmu}).

\medskip
(b) The estimate (\ref{estgamA}) follows
from (\ref{PLinfty}) (that is, $\|\cP u\|_{\infty}\leq \|u\|_{\infty}$)
and  (\ref{PL1wmu}) by means of
the Marcinkiewicz theorem
in a quantitative form with effective constant
\cite[Th.~1.3.2]{Grafakos}.
\eop

\medskip
The result (but not proof) of part (b) of Theorem~\ref{pcurve} carries over {\em mutatis mutandis}\
%
%
if the Poisson integral is replaced by the Cauchy integral
\beq{Cauchyint}
 \St u(z)=\int_{\RR} S(x-t,y)\,u(t)\,dt.
\eeq
Note that $(z-t)^{-1}\in L^{p'}$ as a function of $t$ if $y>0$ and $p'>1$. Therefore by the H\"older inequality 
$\St u(z)$ is pointwise defined in $\HH$ if $u\in L^p$, $1\leq p<\infty$.   

\begin{thm}\label{ccurve}
Under the assumptions (\ref{muwpclass}) and $1<p<\infty$, the operator $\St$ 
is a bounded operator from $L^p(\RR)$ to $L^p(d\mu)$ and the uniform norm bound holds: 
\beq{estgamS}
\ba{ll}\dst
\KScurve(p,k_x,k_y)&
\dst
:= \sup_{\mu\in\Wp(\HH,k_x,k_y)}\,\|\St\|_{L^p(\RR)\to L^p(d\mu)} 
\\[3ex]&\dst
\leq
 \KPcurve(p,k_x,k_y)\left(1+\tan\left(\f{\pi}{2}\,\max(p^{-1},p'^{-1})\right)\right).
\ea
\eeq
\end{thm}

\pf
If $u(x)$ is a real-valued function, then the real part of the analytic function $\St u(z)$ is $\cP u$ and the imaginary part is $\cP v$, where $v$ is the {\em conjugate function}\ of $u$
obtained by the {\em Hilbert transform}
\beq{Ht}
 v(x)=\Ht u(x)=p.v.\f{1}{\pi}\int_\RR \f{f(t)}{x-t}\,dt, \quad x\in\RR.
\eeq
(For details see e.g. \cite[\S~4.1.2]{Grafakos}.)\@ To summarize, we have $\St u=\cP\circ(I+i\Ht)u$, and in this  relation the requirement that $u$ be real-valued can of course be dropped due to linearity over $\CC$. 

By a result of Pichorides,  $\|\Ht\|_{L^p(\RR)}=\tan\left(\f{\pi}{2}\,\max(p^{-1},p'^{-1})\right)$
\cite[Remark 4.1.8]{Grafakos}.
Now (\ref{estgamS}) follows from (\ref{estgamA}).
\eop

\medskip
In the particular case $p=2$ the estimate (\ref{estgamS}) yields 
\beq{estgamS2}
\KScurve(2,k_x,k_y)\leq 4\sqrt{6k_x+k_y}.
\eeq
We allowed a generous simplification in (\ref{estgamA}): $2\pi^{-1}\mapsto 1$.

\subsection{Hausdorff-Young type inequality 
with well-projected weights}
\label{sec:uniflapcurve} 

We know from Theorem~\ref{lrayp} that the Laplace transform
$$ u(t)\mapsto \Lt u(z)=\int_0^\infty u(t) e^{-tz}\,dt, \quad \Re
 z\geq 0,
$$
is a bounded operator from 
$L^p(\RR_+)$ to the space
$L^{p'}(\Ray{\theta})$, where $1\leq p\leq 2$. 
Here $\theta\in[-\pi/2,\pi/2]$ and the uniform bound (\ref{Ltnormp}) holds for the norm.

In this section we generalize this result, replacing the arclength measure supported on the ray $\Ray{\theta}$
with a measure supported on a sector in the right half-plane $\Cp$ and well-projected relative to the sector's bounding rays. 


\newcommand{\kalphax}{k}
\newcommand{\kalphay}{k'}
\newcommand{\kbetax}{l}
\newcommand{\kbetay}{l'}

Let $-\pi/2\leq\al<\be\leq \pi/2$.
Suppose $\mu$ is a measure supported on the sector
$\al\leq\arg z\leq\be$ and well projected relative to its bounding rays (see Definition~\ref{def:alpha-wp}, p.~\pageref{def:alpha-wp}):
\beq{wpmessect}
\ba{l}
\dst \mu\in \Wp_\al(\kalphax,\kalphay),\\[1ex]
\dst \mu\in \Wp_\be(\kbetax,\kbetay).
\ea
\eeq 
The class of all such measures will be denoted
$\Wp_{\al,\be}(\kalphax,\kalphay, \kbetax,\kbetay)$.

\begin{thm}
\label{lmasterthm}
For any $\mu \in \Wp_{\al,\be}(\kalphax,\kalphay, \kbetax,\kbetay)$ and any
$1\leq p\leq 2$, the Laplace transform 
is a bounded operator from $L^p(\RR_+)$ to $L^{p'}(d\mu)$. Moreover, the norm bound is uniform for all measures with given projection constants:   
\beq{Lpmaster}
\ba{l}\dst
\KLcurve(p; \kalphax,\kalphay, \kbetax,\kbetay):=
\sup_{-\f{\pi}2\leq\al<\be\leq \f{\pi}2}\;\;
\sup_{\mu\in\Wp_{\al,\be}(\kalphax,\kalphay, \kbetax,\kbetay)}
\|\Lt\|_{L^p(\Rplus)\to L^{p'}(d\mu)}
\\[4ex]
\dst
\qquad\qquad
\leq (8\pi)^{1/p'}\,
\left(\sqrt{6\kalphax+\kalphay}+\sqrt{6\kbetax+\kbetay}\right)^{2/p'}.
\ea
\eeq
\end{thm}

%
\pf
The inequality (\ref{Lpmaster}) follows by the Riesz-Thorin theorem, cf.\ (\ref{K2interp}), from its two extreme cases:
the $L^1\to L^\infty$ estimate $\KLcurve(1,\dots)\leq 1$, which is trivial since 
$\,|\Lt f(z)|\leq \|f\|_{L^1(\RR_+)}$, and a nontrivial $L^2(\Rplus)\to L^2(d\mu)$ estimate, yet to be established,
\beq{Lcurve2}
\|\Lt f\|_{L^2(d\mu)}\leq \KLcurve(2; \kalphax,\kalphay, \kbetax,\kbetay)\,\|f\|_{L^2(\RR_+)}
\eeq
with constant
\beq{Lgnum2}
\KLcurve(2; \kalphax,\kalphay, \kbetax,\kbetay)\leq  (8\pi)^{1/2}\,
\left(\sqrt{6\kalphax+\kalphay}+\sqrt{6\kbetax+\kbetay}\right).
\eeq
In the proof that follows we could have kept $p$ as a parameter throughout; the only simplification gained by setting $p=2$ is the visual and typographical advantage of the shorter constant in (\ref{estgamS2}) compared to that in (\ref{estgamS}). 
The proof would yield $L^p(\Rplus)\to L^{p'}(d\mu)$ boundedness of the operator $\Lt$ for all $p\in(1,2]$ without resorting to interpolation in the end. 
However, in the corresponding generalization of  (\ref{LcurveKp}) we would end up with a constant that grows as $(1-p)^{-2}$ when $p\to 1^+$, cf.~(\ref{estgamS}).
On the contrary, the right-hand side of (\ref{Lpmaster}) is uniformly bounded for $p\in [1,2]$.
This is why the interpolatory approach is superior when $p$ is close to 1.

\smallskip
We will derive an $L^2(\Rplus)\to L^{2}(d\mu)$ estimate for the operator $\Lt$ from Theorem~\ref{lrayp} using analyticity of the function $\Lt f(z)$ and Theorem~\ref{ccurve}.

\begin{lem}
\label{csector}
Suppose $f(t)$ is a simple function (compactly supported, piece-wise constant, with finite number of values) defined on $\Rplus$. Let 
$-\pi/2\leq\al<\arg z<\be\leq\pi/2$.
We have the 
Cauchy integral representation 
\beq{Cauchy}
\Lt f(z)=\f{1}{2\pi i}\left(
\int_{0}^{\infty} \frac{(\Lt_{\be}f)(x)}{xe^{i\be}-z}\,e^{i\be}\,dx
-
\int_{0}^{\infty} \frac{(\Lt_{\al}f)(x)}{xe^{i\al}-z}\,e^{i\al}\,dx
\right).
\eeq
\end{lem}

\pf
The function $f$ can be written in the form
$$
 f(t)=\sum_{n=1}^N c_n \chi_{[0,b_n]}(t),
$$
where $\chi_S$ denotes the characteristic function of the set $S$. 
The Laplace transform of $f(t)$ is
$$
 \Lt f(\zeta)=\sum_{n=1}^N c_n \zeta^{-1}\,\left(1-e^{-\zeta b_n}\right).
$$
Thus, $\Lt f(\zeta)\sim C\zeta^{-1}$ as $|\zeta|\to \infty$ uniformly for $\al\leq \arg \zeta \leq \be$.
Hence the Cauchy integral along the arc
$\Ga_R=\{\zeta:\,|\zeta|=R, \;\al\leq \arg \zeta\leq \be\}\,$ is estimated as
$$
\left|\frac{1}{2\pi i}\int_{\Ga_R} \frac{\Lt f(\zeta)}{\zeta-z}\,d\zeta\right|
=O(R^{-1}), \quad R\to\infty.
$$
Applying the Cauchy integral formula with
contour of integration consisting of the two straight segments $\{\zeta: 0\leq|\zeta|\leq R,\;\arg\zeta=\al\;\;{\rm or}\;\;\be\}$
 and the arc $\Ga_R$,
we obtain (\ref{Cauchy}) in the limit $R\to\infty$.
\eop


\medskip
Continuing the proof of Theorem~\ref{lmasterthm}, we
estimate the two integrals in the right-hand side of~(\ref{Cauchy}). To be specific, pick the first integral (with $\beta$). 
By Theorem~\ref{lray2}, see~(\ref{L2norm}), we have 
$$
\|\Lt_{\be}f\|_{2}\leq (2\pi)^{1/2}\|f\|_2.
$$
Rotate the coordinate system so as to make the ray of integration $\Ray{\be}$ the positive
$x$ axis and extend the function $\Lt_{\be}f(x)$ by zero for $x<0$.
We are now in the situation of Theorem~\ref{ccurve}, so 
$$
\left\|\f{1}{\pi i}\int_0^\infty\f{(\Lt_\be f)(x)}{xe^{i\be}-z}\,e^{i\be}\,dx\right\|_{L^{2}(d\mu(z))}\leq 4\sqrt{6\kbetax+\kbetay}
\; (2\pi)^{1/2}\|f\|_2.
$$
The estimate for the $\al$-part of the right-hand side of (\ref{Cauchy}) is completely analogous and we conclude 
\beq{LcurveKp}
 \|\Lt f\|_{L^2(d\mu)}\leq 2 (2\pi)^{1/2}\,
\left(\sqrt{6\kalphax+\kalphay}+\sqrt{6\kbetax+\kbetay}\right)\,
\|f\|_2,
\eeq
which coincides with (\ref{Lcurve2}), (\ref{Lgnum2}).

To obtain the estimate (\ref{LcurveKp}), we have so far assumed that $f$ is a simple  function (see Lemma~\ref{csector}).
This constraint is unimportant due to density of simple functions in $\Lp$. 
%
%
The proof of Theorem~\ref{lmasterthm} is complete.
\eop
  
\bigskip\noindent
{\bf Remark}.
In our earlier paper \cite[Note 2 after Theorem 2.1]{SM-CMA10}, a counterpart of Theorems \ref{lrayp} (with an unimportant restriction $p=2$ and an unspecified constant) was found as a particular case of the main result,
which, in terms of logical dependencies,
corresponds to Theorem~\ref{lmasterthm} here. To the contrary, in the present paper a proof of the latter statement depends on the former. 
This logical structure allows us, first, to attend to the question of best constants in Theorem \ref{lrayp}, and, second, to use the Cauchy integral representation with integration along rays (Lemma \ref{csector}), which is 
easier to justify than a similar representation in the Paley-Wiener theory \cite[formula (3.1)]{SM-CMA10}, where the integration path is the imaginary axis.


\section{Applications} 

\label{sec:app}
\newsectnum

\subsection{Uniform Hausdorff-Young estimates in some classes of curves}
\label{sec:unifclasses}

The spaces $L^p(\ga)$ in this section are defined with respect to the arclength measure $\mu_\ga$ on the ordinary rectifiable or compound  curve $\ga$.

\medskip
We begin with a theorem which is the most straightforward extension of the result of our paper \cite{SM-CMA10}.

\medskip
Let $\MonCP$ be the class of monotone curves in $\Cp$ defined in Example~\ref{ex:moncurve} of Sect.~\ref{sec:wpexamples}.
Let $\ConvCP$ be the class of convex curves in $\Cp$ defined in Example~\ref{ex:convcurve}.

\begin{thm}
\label{thm:unifcurv}
If $\ga$ is a monotone or convex curve in $\Cp$ and $1\leq p\leq 2$, then the operator $\Lg$ acts continuously from $L^p(\Rplus)$ to $L^{p'}(\ga)$.
Moreover, we have the uniform norm estimates
\beq{unifmon}
\Kmon(p):=\sup_{\ga\in\MonCP}\,\|\Lg\|_{L^p(\Rplus)\to L^{p'}(\ga)}
\leq (224\,\pi)^{1/p'} 
\eeq
and
\beq{unifconv}
\Kconv(p):=\sup_{\ga\in\ConvCP}\,\|\Lg\|_{L^p(\Rplus)\to L^{p'}(\ga)}
\leq (16\cdot 224\,\pi)^{1/p'} . 
\eeq
\end{thm}

\pf
This theorem is an immediate corollary of Theorem~\ref{lmasterthm}.
If $\ga$ is a monotone curve, then 
$\MonCP\subset\Wp_{-\pi/2,\pi/2}(1,1,1,1)$ by Example~\ref{ex:moncurve}.
Here is the arithmetic connecting (\ref{Lpmaster}) and (\ref{unifmon}):
$$
  8\cdot (2\sqrt{6+1})^2=224.
$$ 
Similarly, if $\ga$ is convex, then 
$\ConvCP\subset\Wp_{-\pi/2,\pi/2}(4,4,4,4)$
by Example~\ref{ex:convcurve},
and (\ref{unifconv}) follows.
\eop

\bigskip
As the next application of Theorem~\ref{lmasterthm} we
obtain uniform estimates for $\Lg$ where $\ga$ is the graph of a Lipschitz function. This result covers some ``wavy'' curves (sinusoids and such), to which
Theorem~\ref{thm:unifcurv} is not applicable. 

\begin{thm}
\label{thm:lipcurve}
(a) Suppose $f(x)$ is a Lipschitz function on an interval $I\subset\Rplus$ with Lipschitz constant $\lam$, see~(\ref{Lipk}).
Let $\ga=\{x+if(x)\mid x>0\}$ be the graph of $f$.
The operator $\Lg$ acts continuously from $L^p(\Rplus)$ to $L^{p'}(\ga)$.
Moreover, the following uniform estimate holds:
\beq{unifgraph} 
\Kgraph(p,\lam):=\sup_{f\;\; \mbox{\footnotesize\rm satisfies (\ref{Lipk})} }\,
\|\Lg\|_{L^p(\Rplus)\to L^{p'}(\ga)}
\leq \big(32\,\pi\,(1+\lam^2)\big)^{1/p'}.
\eeq
(b) If $f\geq 0$ is a Lipschitz function satisfying (\ref{Lipk}) on an interval
$I\subset\RR$ and $\hat\ga=\{f(x)+ix\mid x>0\}$ is the 
graph of $f$ rotated by $90^\circ$, then we have the uniform estimate
\beq{unifgraphy} 
\ba{l}
\dst
\Kgraphy(p,\lam):=
\sup_{\ba{c} \scriptstyle f\geq 0, 
\;\; \mbox{\footnotesize\rm satisfies (\ref{Lipk})}
\ea 
}\,
\|\Lt_{\hat\ga}\|_{L^p(\Rplus)\to L^{p'}(\hat\ga)}
\\[4ex]
\qquad\leq\dst
\big(36\cdot 32\,\pi\,(1+\lam^2)\big)^{1/p'}.
\ea
\eeq
\end{thm} 

\pf
According to Example~\ref{ex:lipcurve}, 
$\ga\in\Wp(\sqrt{\lam^2+1},0)$. Therefore,
Theorem~\ref{lmasterthm} is applicable with
$\al=-\pi/2$, $\be=\pi/2$,
$\kalphax=\kbetax=0$, and $\kalphay=\kbetay=\sqrt{\lam^2+1}$.
Thus we come to (\ref{unifgraph}).

Interchanging the $x$ and $y$ axes, we acquire the extra factor
$6^{2/p'}$ in (\ref{unifgraphy}), which comes from (\ref{Lpmaster}). 
\eop

\medskip
Theorem~\ref{thm:lipcurve} covers some oscillatory curves with unbounded amplitude --- those for which, roughly speaking, the amplitude grows not faster than the period. Example: 
let $0<\al<1$ and 
\beq{oscilcurve}
y=x^\al\,\cos x^{1-\al},\quad x\geq 1.
\eeq 
However, if we replace the arcs of the curve 
(\ref{oscilcurve}) between its
zeros by semicircles (see Fig.~\ref{fig:semicircles}
for $\al=0$),
Theorem~\ref{thm:lipcurve} does not help because the Lipschitz condition no longer holds.
Our next theorem is aimed at curves like that. 
   
Recall the definition of boxed curves and the classes
$\mathfrak{B}(\nu,c)$, see Example~\ref{ex:boxedcurve}
and (\ref{moncomb})--(\ref{hcombc}).
(One can consider mirrored boxes, where $-b_n\leq\Im z\leq 0$, in the same manner).

Introduce subclasses 
$\mathfrak{B}_{\al}(\nu,c)\subset\mathfrak{B}(\nu,c)$.
Here $0<\al<\pi/2$ and we require a curve  of class $\mathfrak{B}_{\al}(\nu,c)$ to lie in the sector
$0\leq\arg z\leq\al$.

To avoid cumbersome expressions,
in the last part of Theorem~\ref{thm:boxcurve} we take an arbitrary but representative particular case: $\al=\pi/4$, 
$\nu=1$. Denote the class of such curves
$$
\mathfrak{B}_*(c)=\mathfrak{B}_{\pi/4}(1,c).
$$

\begin{thm}
\label{thm:boxcurve}
Given $0<\al<\pi/2$, a positive integer $\nu$,
and real $c>0$, let $\ga$ be a boxed curve, $\ga\in\mathfrak{B}_{\al}(\nu,c)$. 
The operator $\Lg$ acts continuously from $L^p(\Rplus)$ to $L^{p'}(\ga)$. Its norm is bounded by a constant that
depends on $\al$, $c$, $\nu$, and $p$.
In particular, the following uniform estimate holds:
\beq{unifbox} 
\Kbox(p,c):=\sup_{f\in\mathfrak{B}_*(c)}\,
\|\Lg\|_{L^p(\Rplus)\to L^{p'}(\ga)}
<
\big(1344\,\pi\,(c^{-1}+2)\big)^{1/p'}.
\eeq
\end{thm}

\pf This theorem follows directly from Theorem~\ref{lmasterthm} and Lemma~\ref{lem:boxproj}.
To prove (\ref{unifbox}),
take $\al=\pi/4$ and $\nu=1$ in (\ref{wpboxk}),
and find the projection constants
$\kappa_1=\kappa_2=4(c^{-1}\sqrt{2}+3)<6(c^{-1}+2)$.
The common value of $\kappa_1$ and $\kappa_2$
is then used as the common value of the constants
$\kalphax,\kalphay,\kbetax,\kbetay$ in (\ref{Lpmaster}).
Substitution yields 
$$
(\sqrt{\quad}+\sqrt{\quad})^2
\leq 4\cdot 7\, \kappa_1 < 168(c^{-1}+2),
$$
and we obtain (\ref{unifbox}). 
\eop

\begin{figure}
\centerline{
\begin{picture}(150,230)
\put(0,0){\includegraphics{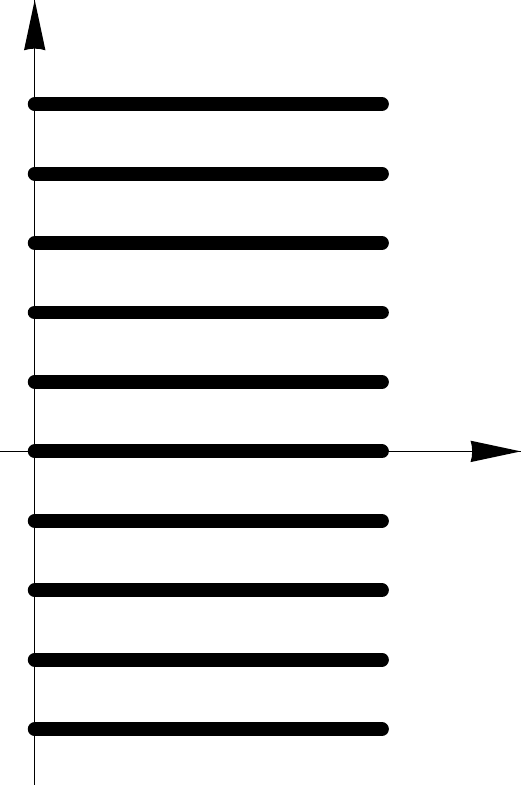}}
\put(139,104){$x$}
\put(20,220){$y$}
\end{picture}
}
\caption{An equispaced vertical comb (Theorem~\ref{thm:vertcomb})}
\label{fig:vertcomb}
\end{figure}

\medskip
Horizontal combs (Example~\ref{ex:horcomb}), being a particular case of boxed curves fall into the scope of 
Theorem~\ref{thm:boxcurve} (under the appropriate assumptions about the parameters $\{a_n\},\{b_n\}$). 
A suitable modification of the theorem, which we omit, handles combs and more general boxed curves rotated through an angle $\al\in(0,\pi/2)$ clockwise or counterclockwise.
However, no such modification is possible for $\al=\pm\pi/2$.
In particular, the equispaced 
vertical comb (Fig.~\ref{fig:vertcomb}) remains out of reach for corollaries of Theorem~\ref{lmasterthm}. 
Well-projected estimates in this case are available for
non-vertical and non-horizontal directions, but we don't have an enclosing sector with suitable boundaries.

Using an extension of technique that led to Theorem~\ref{lmasterthm}, we will prove a Hausdorff-Young type theorem for the 
equispaced
vertical comb. (Theorem~\ref{thm:vertcomb}) in the last section
\ref{sec:vcomb}.

Now let us demonstrate that the relation (\ref{hcombc}) is ``of the right order'', that is, in a certain sense
it is rather close to optimal. 

\subsection{A remark on weighted estimates for the Laplace transform and tightness of conditions in Theorem~\ref{thm:boxcurve}}

Let $w(t)$ be a positive weight function on $\Rplus$
and $1<p\leq 2$.
Bloom's result --- a particular case of Theorem 2(a) in \cite{Bloom} --- gives a necessary and sufficient condition for the boundedness of the Laplace transform as an operator from $L^{p}(\Rplus)$ to
$L^{p'}(\Rplus,w(x)dx)$. The condition is
\beq{Bloomcond}
\Lt w(x)\leq\f{C}{x},\qquad \forall x>0.
\eeq
%
%
In particular, if $w(x)$ is monotone (nonincreasing or nondecreasing), then, as is easily
seen, $x\,\Lt w(x)=O(1)$ as $x\to 0^+$ if and only if 
$w(x)$ is bounded. Thus, for a monotone weight $w$, the necessary and sufficient condition for the estimate 
$$
 \int_0^\infty |\Lt f(x)|^{p'}\, w(x)\,dx\leq
C\|f\|_p^p
$$
to hold is $w(x)=O(1)$. (The result is expected and is not difficult to prove directly.)

Let us examine this fact in connection with our Theorem~\ref{thm:boxcurve}. Specifically, let us consider a horizontal comb 
$\ga$
with teeth $\ga_n=[a_n, a_n+ib_n]$, $n=1,2,\dots$.
Suppose, postponing any effort to justify, that
the absolute values of $\Lt f$ along the tooth $\ga_n$ are comparable with absolute values of $\Lt f$ along the interval
$(a_{n-1},a_n]$ (we take $a_0=0$):
\beq{vagueapprox}
 |\Lt f(a_n+iy)| \approx \left|\Lt f\left(a_n-\frac{a_{n}-a_{n-1}}{b_n}y\right)\right|.
\eeq
We deliberately leave a precise meaning of the $\approx$ sign undefined.%
\footnote{A wishful interpretation: $U\approx V$,
if $c_1\leq U/V\leq c_2$, where $U$ and $V$ are real positive functions, and $c_1$ and $c_2$ are  positive constants.
But it is too good to be true.} 
In such intuitive sense we obtain
$$
\intl_0^{b_n} |\Lt f(a_n+iy)|^p\,dy \approx
\intl_{a_{n-1}}^{a_n} |\Lt f(x)|^p \,w(x)\,dx,
$$
where for $x\in(a_{n-1},a_n]$ 
\beq{combweight}
 w(x)=w_n=\frac{b_n}{a_{n}-a_{n-1}}.
\eeq
The condition $w(x)=O(1)$ can be informally stated as a rule of thumb:
\begin{quote}
\label{aspectratiorule}
{\em The ratio {\em `tooth height to inter-teeth spacing'}\
should be bounded for the Hausdorff-Young inequalities on a comb to hold}. 
\end{quote}
Compare to the condition (\ref{hcombc}) with $k=1$!

\medskip
For instance, if $a_n=n^\al$, $\al\geq 1$, and
$b_n=n^\be$, then 
$w_n\sim \al n^{\al-\be-1}$, so 
$w(x)\sim c x^{1-(\be+1)/\al}$ for large $x$. 
Thus, by a simple heuristic argument, the Hausdorff-Young theorem for $\Lg$ is expected to be true if and only if $\be\leq \al-1$. The ``if'' statement is contained in Theorem
\ref{thm:boxcurve} and the ``only if'' will be proven by making the vague relation (\ref{vagueapprox}) more mathematical.

A similar heuristic analysis with the same conclusion
can be carried out for more general boxed curves.  


As promised, we will now show that the Hausdorff-Young theorem fails for the comb $\ga$ as above in the case $\be> \al-1$.
It suffices to construct a family of functions $\{f_k\}$,
depending on parameter $k>1$ such that 
\\[1ex]
(a) \hspace{1em} $\|f_k\|_p=1$,
\\[1ex]
(b) \hspace{1em} 
$\int_0^\infty\,|\Lt f_k|^{p'}\,w(x)\,dx\,\to\,\infty\,$
as $k\to\infty$,
and 
\\[1ex]
(c) \hspace{1em} all $f_k$ satisfy a rigorous replacement of (\ref{vagueapprox}) --- the inequality
$$
 |\Lt f_k(a_n+iy)| \geq c_1 \left|\Lt f_k\left(a_n-\frac{a_{n}-a_{n-1}}{b_n}y\right)\right|,
\quad 0\leq y\leq b_n
$$
with some $c_1>0$ independent of $n$ and $k$ .

Let $\delta=1-(\be+1)/\al\,>0$. Fix $p\in(1,2]$.
Define $f_k(x)=A_k e^{-kx}$ with $A_k=(kp)^{1/p}$. The choice of $A_k$ is dictated by condition (a). Next, $\Lt f_k(t)=A_k (t+k)^{-1}$, so
$$
\ba{l}
\dst
 \int_0^\infty \left|\Lt f_k(t)\right|^{p'}\,(1+t)^\delta\,dt
\;>\;
(kp)^{p'/p}\,\int_0^\infty \f{t^\delta\,dt}{(t+k)^{p'}}
\\[4ex]\dst
\quad=\;k^{p'/p+\delta+1-p'}\,p^{p'/p}\,\int_0^\infty
\f{\tau^\delta\,d\tau}{(\tau+1)^{p'}}
\;=\;\; C k^\delta,
\ea
$$
where $C$ does not depend on $k$. Thus condition (b) holds true. 

To simplify analysis of condition (c),
let $t=a_n-(a_{n}-a_{n-1})y/b_n$ and $t+\Delta t=a_n+iy$.
Notice that $|\Delta t|\leq ct$, where $c$ depends only on $\alpha$. (A stronger estimate $|\Delta t|\leq ct/n$ is true but unnecessary.)
The condition (c) holds with constant $c_1=(1+c)^{-1}$.
Indeed,
$$
 \f{|\Lt f_k(t)|}{|\Lt f_k(t+\delta t)|}
=\left|\f{t+\Delta t+k}{t+k}\right|
\leq 1+\f{\Delta t}{t+k}
< 1+\f{\Delta t}{t}\leq 1+c.
$$

Conclusion: A comb with parameters $a_n=n^\al$, $b_n=n^\be$ satisfies conditions of Theorem \ref{thm:boxcurve} (assuming $\al\geq 1$) if and only if $\beta\leq\alpha-1$. The boundedness of the operator $\Lg$ from $L^{p}(\Rplus)$ to
$L^{p'}(\ga)$ takes place in the same range of values of $\be$. This example should help to get an intuitive feeling about tightness of conditions in Theorem~\ref{thm:boxcurve}.

\subsection{Two maximal theorems}
\label{sec:maxthms}

\unl{Maximal Paley-Wiener theorem }

\medskip
Let $f$ be a function on $\Rplus$ for which the Laplace
transform $\Lt f$ is defined in the right half-plane $\Cp$. 
Define the {\em maximal Laplace transform}\ of $f$,
\beq{maxLt}
 \Lt^* f(y)=\sup_{x>0} |\Lt f(x+iy)|.
\eeq

\begin{thm}
\label{lmaxPW}
For $1<p\leq 2$ the inequality 
\beq{maxPW}
 \left(\int_{-\infty}^\infty |\Lt^* f(y)|^{p'}\,dy \right)^{1/p'}
\leq \KmaxPW(p)\, \|f\|_p
\eeq
holds. (In particular, $|\Lt^* f|<\infty$ a.e.) \@ Here
$$
 \KmaxPW(p)\leq (192\pi)^{1/p'}.
$$
\end{thm}

\noindent
{\bf Remark}. If instead of $\Lt^* f(y)$ we consider the function $\Lt^{**} f(x)=\sup_{y} |\Lt f(x+iy)|$, then there is an obvious pointwise  estimate
$$
 |\Lt^{**} f(x)|\leq (\Lt|f|)(x),
$$
and the corresponding ``maximal theorem'' is a trivial corollary of the Hausdorff-Young inequality for the Laplace transform on $\Rplus$. 

\begin{lem} 
\label{maxPW2}
Theorem \ref{lmaxPW} allows the following equivalent formulation:
\\ Let $1<p\leq 2$. 
The inequality
\beq{PWb}
 \left(\int_
{-\infty}^\infty 
|\Lt f(b(y)+iy)|^{p'}\,dy \right)^{1/p'}
\leq \KmaxPW(p) \|f\|_p
\eeq
holds for any $f\in L^p(\Rplus)$ and any simple (i.e.\ piece-wise constant, finitely-valued) function $b:\,
\RR\to\Rplus$.
\end{lem}

\pf Since the left-hand side of (\ref{PWb}) is not greater than the l.h.s.\ of (\ref{maxPW}), the original formulation implies the new one. 
Let us prove the converse. 

Suppose 
$b:\,\RR
\to\Rplus$ is a measuruable function. 
Approximate $b(y)$ by simple functions
$b^{(n)}(y)\to b(y)$ a.e. 
For each of the functions $b^{(n)}$ we have the inequality (\ref{PWb}).
Since $\Lt f(z)$ is continuous, $\Lt f(b^{(n)}(y)+iy)\to \Lt f(b(y)+iy)$ a.e. 
%
%
By Fatou's Lemma, the inequality (\ref{PWb}) still holds for $b(y)$.
Thus  
the class of admissible functions $b$ 
in Lemma~\ref{maxPW2}
has been extended from simple to measurable. 

\smallskip
Next, fix $f\in L^p(\Rplus)$. The function $\Lt f(z)$ is analytic in $\Cp$ and for any $y$,
$|\Lt f(x+iy)|\to 0$ as $x\to+\infty$. There are well-defined functions $\Lt_{\eps}^* f(y)=\max_{x\geq\eps} |\Lt f(x+iy)|\;$ 
and 
$$
b_\eps(y)=\min\left\{x\geq\eps\;|\;\Lt_{\eps}^* f(y)=|\Lt f(x+iy)|\,\right\}.
$$
The function $b_\eps(y)$ (which depends on $f$) is 
obviously
measurable. 
By the previous, the estimate (\ref{PWb}) with $b_\eps$ in place of $b$ holds true. The right-hand side is independent of $\eps$. 
Finally, $\Lt_{\eps}^* f(y)\to \Lt^* f(y)$ as $\eps\to+0$ and the inequality (\ref{maxPW}) follows 
by the Monotone Convergence Theorem.
\eop 

\begin{figure}
\centerline{
\begin{picture}(170,270)
\put(0,0){\includegraphics{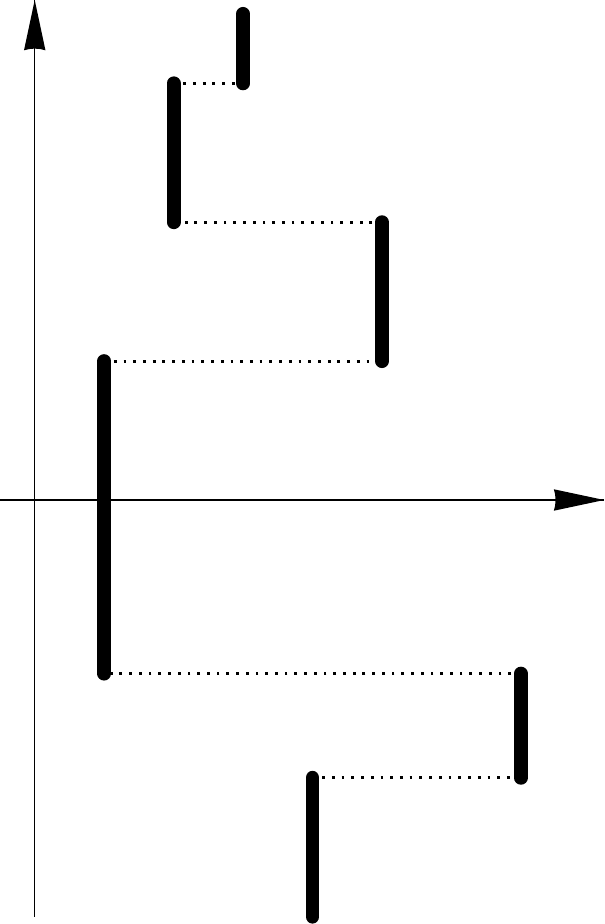}}
\put(164,130){$x$}
\put(20,258){$y$}
\end{picture}
}
\caption{A simple discontinuos vertical curve for the Maximal Paley-Wiener Theorem (Th.~\ref{lmaxPW})}
\label{fig:maxPW}
\end{figure}

\medskip
\noindent{\it Proof of Theorem~\ref{lmaxPW}}. We 
prove the Theorem in the form 
stated in Lemma~\ref{maxPW2}.
Consider a simple positive function $b(y)$ defined for almost all $y\in\RR$,
$$
 b(y)=b_n\qquad \mbox{\rm for}\quad a_{n-1}< y< a_{n},\quad n=1,\dots, N,
$$
where $a_0=-\infty$, $a_N=+\infty$.
To it, there corresponds a compound curve  $\ga=\{\ga_n\}$,
where $\ga_n$ is the vertical interval 
$(b_n+ia_{n-1},\; b_n+ia_n)\;$ (Figure~\ref{fig:maxPW}).

It remains to apply Theorem~\ref{lmasterthm} to the
arclength measure $\mu_\ga$=$|dy|_{\ga}$. 
Clearly,
$\mu_\ga\in \Wp_{-\pi/2,\pi/2}(1,0,1,0)$
--- cf.~(\ref{wpmessect}).
We 
get 
(\ref{PWb}) with 
constant
$$
 \KmaxPW(p)\leq K_5(p,1,0,1,0)\leq (8\cdot 2^2\cdot 6\pi)^{1/p'},
$$
as claimed. 
\eop

\bigskip
\unl{A restricted Maximal Angular Theorem}

\medskip
The result of this section, Theorem~\ref{restrangmax}, viewed separately, may look very artificial. It should be considered in an appropriate context. The context is what we call the ``Maximal Angular Theorem'' (MAT).

In this paper we prefer to define a status of MAT not in terms of its claimed or conjectured trueness, but in terms of goal setting --- as a ``desirable result''.  This will be signified by the question mark over inequality sign.  
We hope to clarify the status and to discuss further details in a separate paper. 

Let us explain {\em what}\ is desired.

By analogy with (\ref{maxLt}) 
introduce the {\em angular maximal Laplace transform} 
$$
 (\Lt^\sharp f)(r)=\sup_{|\theta|<\pi/2} |\Lt f(re^{i\theta})| 
,\quad r>0 .
$$

The most straightforward form of MAT would assert the Maximal Angular Inequality. It is but the Hausdorff-Young type inequality for $\Lt^\sharp f$,
\beq{maxangpure}
 \left(\int_{0}^\infty |\Lt^\sharp f(r)|^{p'}\,dr \right)^{1/p'}
\;
\stackrel{\raisebox{3pt}{?}}{\leq} 
\;C(p)\, \|f\|_p,\qquad 1<p\leq 2.
\eeq 

By analogy with Lemma~\ref{maxPW2}, an
equivalent formulation that
does not involve the operator $\Lt^\sharp$
can be given.
%
Recall the notion of a simple discontinuous radial curve introduced in Example~\ref{ex:maxang} of Section~\ref{sec:wpexamples}.

So, the second formulation of MAT would propose that
{\it for any $p\in(1,2]$ and any simple discontinuous radial curve $\ga$  
in $\Cp$ with polar equation $\theta=\al(r)$, the inequality
\beq{maxang}
\left(\int_0^\infty|\Lt f(re^{i\al(r)})|^{p'}\,dr
\right)^{1/p'}\;
\stackrel{\raisebox{3pt}{\rm ?}}{\leq} 
\;C(p)\,
\|f\|_p
\\[2ex]
\eeq
holds with constant $C(p)$ that does not depend on functions $f$ and $\al$.}

\medskip
Theorem~\ref{lrayp} (apart from the value of 
the constant) is contained in the inequality (\ref{maxang})  
as a very special case: $\al(r)=\const$. The relation between Theorem~\ref{lrayp} and (\ref{maxang}) parallels
the relation between the classical Paley-Wiener theorem 
and our Theorem~\ref{lmaxPW} in the form (\ref{PWb}).

\medskip
The theorem below is a prototype of MAT, with constraints
that may turn out to be unnecessary.

Given a positive integer $\nu$ and a real $c>1$, 
introduce the class $\mathfrak{R}(\nu,c)$ of simple discontinuous radial curves that lie in the right half-plane and satisfy the ``long radial runs'' condition (\ref{rgeom}).
 
Unlike in the definition of the classes $\mathfrak{R}(\ph,\nu,c)$ in Example~\ref{ex:maxang}, the new classes
do not depend on the sector aperture $\ph$.
Clearly, $\mathfrak{R}(\nu,c)\supset \mathfrak{R}(\ph,\nu,c)$
for all $\ph\in(0,\pi/2)$.

\begin{thm}
\label{restrangmax}
Suppose $\nu$ is a positive integer and $c>1$.
For any 
%
$p\in [1,2]$ the supremum
$$
\Kmaxang(p,\nu,c):=\sup_{\ga\in\mathfrak{R}(\nu,c)}
\|\Lg \|_{L^p(\Rplus)\to L^{p'}(\ga)}
$$
taken over all simple discontinuous radial curves
$\ga$ of the class  $\mathfrak{R}(\nu,c)$ 
is finite and
\beq{rmaxang}
\Kmaxang(p,\nu,c)\leq \left(18\cdot 8\cdot 24\,\left(\f{\log 2}{\log c}+\nu+2\right)\right)^{1/p'}.
\eeq
\end{thm}

\pf
If $\ga\in\mathfrak{R}(\pi/3,\nu,c)$, then
according to Lemma~\ref{lem:wprad}, $\ga$ is well-projected, more precisely 
$
\ga\in \Wp(\kappa,0),
$
with projection constant
\beq{kappapi3}
\kappa=2\left(\f{\log 2}{\log c}+\nu+2\right).
\eeq
obtained by substitution of $\sec\pi/3=2$ to 
(\ref{wpradk}).

By symmetry, $\ga$ is also well-projected in the coordinate system rotated through the angle $\pi/3$ counterclockwise:
$\ga\in\Wp_{\pi/3}(\kappa,0)$. 

Thus we have the inclusions (\ref{wpmessect}), where
$\mu=\mu_\ga$ (which is the radial Lebesgue measure $dr$),
$\al=0$, $\be=\pi/3$, $\kalphax=\kbetax=\kappa$, and $\kalphay=\kbetay=0$.
Therefore, 
$$
 \ga\in\Mes_{0,\pi/3}(\kappa,0,\kappa,0).
$$ 

Since the ``long radial runs'' condition (\ref{rgeom})
is rotation-invariant, a similar well-projectedness property 
is true for any sector $\Phi_{\al,\al+\pi/3}=\{z\mid\al\leq\arg z\leq \al+\pi/3\}$ with
$-\pi/2\leq\al\leq\pi/6$.

That is, if $\ga\in\mathfrak{R}(\nu,c)$ and $\ga$ lies in the
sector $\Phi_{\al,\al+\pi/3}$, then 
\beq{secpi3}
 \ga\in\Mes_{\al,\al+\pi/3}(\kappa,0,\kappa,0).
\eeq

Let us split a given curve $\ga\in\mathfrak{R}(\nu,c)$ into three parts 
$$
\ga=\ga_{-1}\cup \ga_0\cup \ga_1,
$$ 
so that
$\ga_j $ lies in the sector $\Phi_{-\pi/6+j\pi/3,\,\pi/6+j\pi/3}$, $j=-1,0,1$.


Taking $\ga=\ga_j$ 
in (\ref{secpi3}) and applying Theorem~\ref{lmasterthm} we obtain
$$
 \|\Lt\|_{L^p(\Rplus)\to L^{p'}(\ga_j)}
\leq 
\KLcurve(p; \kappa, 0, \kappa, 0)
\leq
(8\pi\cdot 4\cdot 6\kappa)^{1/p'}.
$$
Hence for the whole curve $\ga$ we get
$$
 \|\Lt\|_{L^p(\Rplus)\to L^{p'}(\ga)}
\leq 
3\;(8\cdot 24\pi\kappa)^{1/p'}.
$$
A minor improvement --- replacement of the factor 3 by $3^{2/p'}$ --- is possible due to the Riesz-Thorin theorem, 
as before (cf.~(\ref{Lgnum2})$\Rightarrow$(\ref{Lpmaster})).

Substition of the value of $\kappa$ from (\ref{kappapi3}) yields the estimate
(\ref{rmaxang}). 
\eop

\bigskip\noindent
{\bf Remark}.
If $\ga$ is a ray as in Theorem~\ref{lrayp}, then
$\ga\in\mathfrak{R}(1,c)$ with arbitrarily large $c$.
Theorem~\ref{restrangmax} gives a uniform bound for
the constant $K_2(\theta,p)$ in (\ref{Lrp}):
$$
 K_2(\theta,p)\leq 
\limsup_{c\to\infty}
\Kmaxang(p,1,c)\leq 10368^{1/p'}.
$$ 
where $10368=18\cdot 8\cdot 24\cdot 3$.
The constant $2\pi$ in Theorem~\ref{lrayp}, inequality
(\ref{Ltnormp}), is --- not surprisingly --- much sharper.

\subsection{%
Beyond the reach of Theorem~\ref{lmasterthm}:
Vertical comb 
}
\label{sec:vcomb}

A 
vertical comb (defined in Example~\ref{ex:horcomb}) 
is $\al$-well projected for any $\al$ that is not a multiple of $\pi/2$. However, Theorem~\ref{lmasterthm} is not applicable, because there is no suitable enclosing angle --- see (\ref{wpmessect}), the definition of classes $\Mes_{\al,\be}$.  Yet a Hausdorff-Young type theorem 
is valid 
at least for an equispaced
vertical comb 
with teeth of equal length
(Fig.~\ref{fig:vertcomb}). 

\begin{thm}
\label{thm:vertcomb}
Let $\ga$ be the doubly infinite vertical comb with teeth 
$[in, in+b]$, $b>0$, $n\in\ZZ$. For any $p\in[1,2]$
the operator $\Lg$ acts continuously from $L^p(\Rplus)$ to $L^{p'}(\ga)$. That is, the inequality
\beq{Lgvcomb} 
\left(\sum_{n=-\infty}^\infty \int_0^b |\Lt f(x+in)|^{p'}\,dx
\right)^{1/p'}\,\leq\, 
\Kvcomb(b,p) \|f\|_p
\eeq
holds with a constant $\Kvcomb(b,p)$ independent of $f$.
\end{thm} 

We will give two proofs. The first one is direct, short, and   demonstrates rather explicitly how to exploit the oscillatory factors $e^{int}$ --- the only elements in 
the integrand of the Laplace transform that vary among the infinitely many teeth.  
The second proof is much less elementary but allows greater flexibility; it relies on a refinement of Theorem~\ref{pcurve}. 
Further applications of this technique must be left to another publication.

\medskip
{\it First proof}.
Like in Theorem~\ref{lmasterthm}, the $L^p\to L^{p'}$ boundedness of $\Lg$ for $1\leq p\leq 2$ follows by Riesz-Thorin interpolation
from the two extreme cases: $p=1$ (trivial case) and  $p=2$, which is the subject of the rest of this proof.

Given $f\in L^2(\Rplus)$  we must find a bound for
\beq{vcombL2f}
\|\Lg f\|_{L^2(\ga)}^2=\sum_{n=-\infty}^\infty \int_0^b \left|\int_0^\infty e^{-t(x+in)}\,f(t)\,dt\right|^2\,dx.
\eeq
It will require an application of Parceval's identity and the following estimate based on Hilbert's inequality for the norm of the infinite 
matrix with entries $a_{mn}=(m+n+1)^{-1}$, $m,n\geq 0$.

\begin{lem}
\label{lem:gxt}
Suppose $f(t)$ is a finite continuous function on $\Rplus$
and define
$$
 g(x,t)=\sum_{n=0}^\infty f(t+2\pi n)\,e^{-(t+2\pi n)x},
$$ 
where $t\in[0,2\pi)$ and $x\in(0,b)$.
Then 
\beq{gxL2}
\int_0^b\,dx\,\int_0^{2\pi} |g(x,t)|^2\,dt
\leq
C(b)\,\|f\|_2^2
\eeq
with constant $C(b)$ independent of $f$.
\end{lem}

\pf
Note that the series defining $g(x,t)$ terminates after finitely many terms, so a question about convergence does not arise. 
Next, we may assume that $f\geq 0$.
Denote for brevity $f_n(t)=f(t+2\pi n)$. Then
$$
 \int_0^b |g(x,t)|^2\,dx= \sum_{m,n\geq 0} f_m(t)\,f_n(t)\,
\int_0^b e^{-(2t+2\pi(m+n))x}\,dx.
$$
The function $(1+x^{-1})(1-e^{-x})$ is bounded on $\Rplus$. Therefore the inequality 
$$
 \int_0^b e^{-rx}\,dx =\,\f{1-e^{-br}}{r}\,\leq \f{C_1 b}{1+br},
$$
holds for all $r>0$ with some absolute constant $C_1$. 
From this we infer
$$
\int_0^b |g(x,t)|^2\,dx\,
\leq\,
C_2(b)\,\sum_{m,n\geq 0} \f{f_m(t) f_n(t)}{1+m+n}
\;\leq
\;
C_3(b)\sum_{n=1}^\infty |f_n(t)|^2.
$$
The last step is Hilbert's inequality \cite[\S~9.1]{HLP}.
Integrating over $t\in[0,2\pi]$ we come to (\ref{gxL2})
with $C(b)=C_3(b)$.
\eop

\medskip
Using the notation $g(x,t)$ introduced in Lemma \ref{lem:gxt}, write the inner integral in (\ref{vcombL2f})
as
$$
\ba{lcl}
\dst
 \int_0^\infty e^{-t(x+in)}\,f(t)\,dt
&=&\dst
\sum_{m=0}^\infty \int_0^{2\pi} e^{-(t+2\pi m)x-int}
\,f(t+2\pi m)\,dt
\\[3ex]
&=&\dst
\int_0^{2\pi} e^{-int}\, g(x,t)\,dt.
\ea
$$
Substituting this to (\ref{vcombL2f}) and using Parceval's identity we get 
$$
\ba{lcl}
\dst
\|\Lg f\|_{L^2(\ga)}^2&=& \dst
\sum_{n=-\infty}^\infty \int_0^b \left|\int_0^{2\pi} e^{-int}\,g(x,t)\,dt\right|^2\,dx
\\[3ex]
&=&\dst
2\pi\int_0^b dx\int_0^{2\pi} |g(x,t)|^2\,dt.
\ea
$$
Application of Lemma~\ref{lem:gxt}
leads to the inequality 
$$
\|\Lg f\|_{L^2(\ga)}^2\leq 2\pi\,C(b)\,\|f\|_2^2.
$$
The class of admissible functions in Lemma~\ref{lem:gxt}
is dense in $L^2(\Rplus)$, hence the proof is complete. 
\eop

\medskip
{\it  Second proof}.

Theorem~\ref{lmasterthm} was obtained through the chain of dependent statements:  
$$
\mbox{\rm Theorem~\ref{pcurve}(a) }
\;\Longrightarrow\;
\mbox{\rm Theorem~\ref{pcurve}(b) }
\;\Longrightarrow\;
\mbox{\rm Theorem~\ref{ccurve} }
\;\Longrightarrow\;
\mbox{\rm Theorem~\ref{lmasterthm}}.
$$
The proofs of all those implications did not depend on 
well-projectedness. Therefore, if we manage to  
modify proof of Theorem~\ref{pcurve}(a) so as to make it work for the vertical comb, Theorem~\ref{thm:vertcomb} will follow. The next lemma is the desired  accommodation of Theorem~\ref{pcurve}(a). It asserts a weak $L(1,1)$ boundedness of the Poisson transform $\cP$ (\ref{Poissonint}) restricted to a {\it horizontal}\ comb corresponding to our given vertical comb by $90^\circ$ rotation. 

\begin{lem}
\label{pcurve2}
If $v(z)=\cP u(z)$ is the Poisson transform of a function
$u\in L^1(\RR_+)$, then for any $\lam>0$
\beq{PL1wmu2}
\sum_{n=-\infty}^\infty \big|\{y\in(0,b)\,:\; |v(n+iy)|>\lam  \}\big|\;
\leq\;
2(\pi^{-1}+3b) 
\,\f{\|u\|_1}{\lam}.
\eeq
\end{lem}

\pf To understand a necessary modification to 
the proof of Theorem~\ref{pcurve}, let us examine why that proof fails in the present situation.

The estimates (\ref{mesSy}) for $|A_y|$ and (\ref{mesSx}) for $|A_x|$ are valid; they are not bound to any particular curve or well-projected measure. But (\ref{mesSy}) does not help, since the $y$-projection of our comb has infinite multiplicity. And (\ref{mesSx}) plays no role as the comb's $x$-projection has zero length. 

Thus the reason of failure is the unfortunate distribution of weights in (\ref{wpmes}): $k_x=0$ and $k_y=\infty$.

In the modified proof we break $v$ into two parts,
$v_d$ (``diagonal'') and $v_r$ (``residual'').
The restriction $v_d|_{\ga_n}$ will
depend only on the restriction of $u$ on the 
interval $(n-1/2,n+1/2)$, escaping the
infinite multiplicity trouble. 

An estimate for $v_r$ will be due to the fact that
$|v_r(n+iy)|$ can not be big without the
uncentered maximal function of $u$ being correspondingly big
on the whole interval $(n-1/2,n+1/2)$. 

\smallskip
For $h>0$,
define a truncated Poisson kernel
$$
P_h(x,y)=P(\max(h,|x|),y).
$$
Let us write
$$
 v(z)=v_d(z)+v_r(z),
$$
where 
$$
\ba{l}
 v_d(z)=u \;{*}_{x}\, (P-P_{1/2}),
\\[2ex]
 v_r(z)=u \;{*}_{x}\; P_{1/2}.
\ea
$$
We will estimate
the weak norms $\|v_d\|_{L_\ga^{1,\infty}}$ and $\|v_r\|_{L_\ga^{1,\infty}}$
separately.

\medskip
\unl{Estimation of $\|v_d\|_{L_\ga^{1,\infty}}$}.
Let $u_n(x)$ denote the part of the function $u$ belonging to the interval $(n-1/2, n+1/2)$:
$$
 u_n(x)=\left\{\ba{l}
 u(x), \quad x\in \left(n-\f12,\;n+\f12\right),
\\[2ex]
  0,\qquad x\notin \left(n-\f12,\;n+\f12\right).
\ea
\right.
$$ 
Since $(P-P_{h})(x,y)=0$ for $|x|>h$, we have
$$
\ba{l}\dst
 v_d(n+iy)=\int_{-1/2}^{1/2}
(P-P_{1/2})(t,y)\;u(n-t)\,dt
\\[4ex]\dst\qquad\qquad
=\int_{-1/2}^{1/2}
(P-P_{1/2})(t,y)\;u_n(n-t)\,dt.
\ea
$$
Like in the proof of the estimate (\ref{mesSy}), we use the inequalities
$$
 (P-P_{1/2})(t,y) < P(t,y) \leq \f{1}{\pi y}
$$
to infer
$$
\big|\{y\,:\; |v_d(n+iy)|>\lam\}\big|
\leq \f{1}{\pi\lam}\|u_n\|_{1}.
$$ 
Summing over $n$ yields
\beq{ud-L1w}
\|v_d\|_{L_\ga^{1,\infty}}
\leq \f{1}{\pi}\|u\|_{1}.
\eeq

\medskip
\unl{Estimation of $\|v_r\|_{L_\ga^{1,\infty}}$}.
Introduce the sets
\beq{Ynlam}
 Y_n(\lam)=\{y\in (0,b)\mid |v_r(n+iy)|>\lam\}.
\eeq
We have to estimate $\sum_n |Y_n(\lam)|$.
The idea is to reduce the problem of estimation of sizes of 
the ``vertical'' sets $Y_n(\lam)$ to a similar problem for ``horizontal'' sets $D_U(\mu)$ (see~(\ref{Ds})) described by the inequality $U(x)>\mu$ for an auxiliary function ---  the {\em uncentered Hardy-Littlewood maximal function of $u$}, 
$$
 U(x)=\sup_{I\ni x}\f{1}{|I|}\int_I |u(t)|\,dt.
$$
Here $I$ is any finite interval containing $x$.

\smallskip
We will need a sub-lemma.

\begin{lem}
\label{lem:vrU}
There exists a positive decreasing function $F(t)$ defined for $t\geq 0$ such that $F(0)=1$, $F(t)=O(1/t)$ as $t\to\infty$, and for any $x\in (x_0-h,x_0+h)$ and any $y>0$ 
\beq{vr-maxcomp2}
 |(u*P_h)(x_0+iy)|\leq F(h/y)\,U(x). 
\eeq 
\end{lem}

\label{PleqU} 
\pf
Introduce 
$$
 \rho_y(s)=-2s\,\f{\partial P(s,y)}{\partial s}
\,=\,
 \frac{4}{\pi}\,\frac{ys^2}{(y^2+s^2)^2}
$$
and 
$$
F(t)=\int_t^\infty \rho_1(s)\,ds, \quad t\geq 0.
$$
The properties of $F(t)$ asserted in Lemma are easy to verify.

The function $P_h(x,y)$ is a positive linear combination of
the normalized characteristic functions $(2s)^{-1}\chi_{[-s,s]}(x)$, $s\geq h$, taken with weights $\rho_y(s)$:
$$
P_h(x,y)=
\int_{\max(h,|x|)}^\infty \f{\rho_y(s)}{2s}\,ds
=
\int_{h}^\infty \f{\chi_{[-s,s]}(x)}{2s}\,\rho_y(s)\,ds.
$$
The total weight in this representation is 
$\int_h^\infty \rho_y(s)\,ds=F(h/y)$.  

By definition of the maximal function, for $|x-x_0|<s$ we have 
$$
\f{1}{2s}\,\int_{x_0-s}^{x_0+s} \, |u(t)|\,dt\,\leq\,U(x).
$$
Therefore
$$
\ba{rcl}
\dst 
 |(u *P_h)(x_0+iy)|
 & \leq &
\dst
\int_{-\infty}^\infty P_h(x_0-t,y)\, |u(t)|\,dt
\\[3ex]
& = &
\dst
\int_h^\infty \rho_{y}(s)\,ds \,\int_{-\infty}^\infty \f{\chi_{[-s,s]}(x_0-t)}{2s}\, |u(t)|\,dt
\\[3ex]
& = &
\dst
\int_h^\infty \rho_{y}(s)\,ds \;\;
\f{1}{2s}\,\int_{x_0-s}^{x_0+s} \, |u(t)|\,dt
\\[3ex]
& \leq &
\dst
\int_h^\infty \rho_{y}(s)\, U(x)\,ds \,
=\,
F\left(
h/y
\right)\,U(x).
\ea
$$
Q.E.D. 
\eop

\medskip
Returning to proof of Lemma~\ref{pcurve2}, we apply Lemma~\ref{lem:vrU} with $h=1/2$. 
Suppose that $y\in Y_n(\lam)$. 
According to Lemma~\ref{lem:vrU}, if $|x-n|<1/2$, then
\beq{vr2U}
 \lam<|v_r(n+iy)|\leq 
F(1/(2y))\,U(x)\leq F(1/(2b))\,U(x).
\eeq
 
The above inequality is a key to the proof. The situation is indeed quite surprising. Having assumed that $|v_r(n+iy)|>\lam$ for {\em some}\ $y\in (0,b)$, 
we have obtained a lower bound for the maximal function $U(x)$ not just at a single point but for {\em all}\ $x$ in an interval of length 1:
$$
U(x)>\mu=
\f{\lam}{F(1/(2b))},
\qquad
\forall x\in
(n-1/2, \,n+1/2). 
$$
Thus we formalize the aforementioned connection between
the ``vertical'' sets $Y_n(\lam)$ and the ``horizontal'' sets
$\Ds_U(\mu)$:
$$
Y_n(\lam)\neq\emptyset
\qquad\Longrightarrow
\qquad 
(n-1/2, \,n+1/2)\subset \Ds_{U}(\mu).
$$  
The intervals corresponding to different values of $n$ are disjoint, while the set $\Ds_{U}(\mu)$ is the same for all $n$. Therefore
$$
 \bigcup_{n:\;Y_n(\lam)\neq\emptyset}
(n-1/2, \,n+1/2)\subset \Ds_{U}(\mu).
$$
Passing from sets to their Lebesgue measures we get (recall notation (\ref{Df}))
$$
 \#\{n\mid Y_n(\lam)\neq\emptyset\} \leq d_U(\mu).
$$
On the other hand, the (uncentered) Hardy-Littlewood maximal theorem states that \cite[Theorem~2.1.6]{Grafakos}
$$
 d_U(\mu)
\leq 3\f{\|u\|_1}{\mu}.
$$
Consequently,
$$
\#\{n\mid Y_n(\lam)\neq\emptyset\} \leq 
3 
\,\f{\|u\|_1}{\lam}.
$$
Since the length of every $Y_n(\lambda)$ is at most $b$, we 
conclude that
$$
 \sum_{n=-\infty}^\infty |Y_n(\lam)|
\leq b\,\cdot\,\#\{n\mid Y_n(\lam)\neq\emptyset\}
\leq 3b 
\,\f{\|u\|_1}{\mu}
=C(b)\,\f{\|u\|_1}{\lam}
$$
with $C(b)=3b\cdot F(1/(2b))$.
By definition of the weak $1$-norm it follows that  
\beq{ur-L1w}
\|v_r\|_{L_\ga^{1,\infty}}\leq C(b)\,\|u\|_1.
\eeq

To present the end result in a simple form, we now
coarsen the constant: $F(1/(2b))\leq 1$, hence $C(b)\leq 3b$.
We could have omitted the constant $F(1/(2b))$ already in the r.h.s.\ of (\ref{vr2U}) but we preferred to keep it ---  
mostly in order to justify a notational distinction between $\lam$ (the parameter in the ``vertical'' inequality) 
and $\mu$ (the parameter in the ``horizontal'' inequality).

\smallskip
Combining (\ref{ud-L1w}) and (\ref{ur-L1w}), we obtain
(\ref{PL1wmu2}). 
(Recall that $\|\cdot\|_{1,\infty}$ satisfies the quasi-triangle inequality $\|f+g\|_{1,\infty}\leq 2(\|f\|_{1,\infty}+\|g\|_{1,\infty})$, hence the extra factor 2).

\smallskip
The proof of Lemma~\ref{pcurve2}, and with it --- the second proof of Theorem~\ref{thm:vertcomb}, is finished.
\eop

\section*{Acknowledgements}

A.M.\ acknowledges a financial support by CONACYT.
S.S.\ acknowledges a financial support by NSERC.

\pagebreak

\clearpage

\clearpage

\begin{figure}
\centerline{
\begin{picture}(175,300)
\put(0, 0){\includegraphics{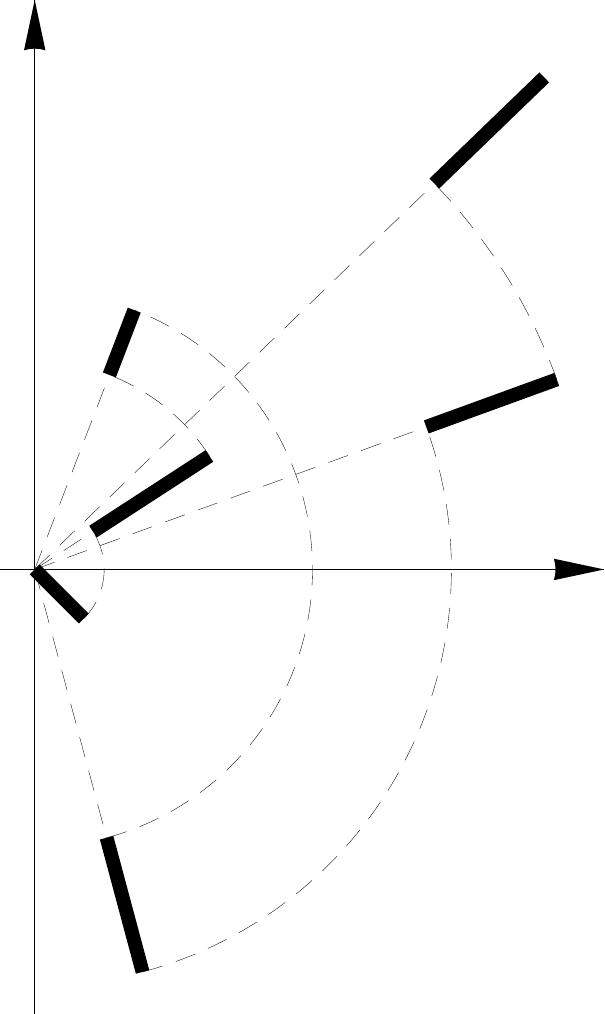}}
\put(164,138){$x$}
\put(20,283){$y$}
\end{picture}
}
\caption{A simple discontinuous radial curve 
(Example~\ref{ex:maxang})}
\label{fig:maxang}
\end{figure}
\clearpage

\clearpage

\clearpage

\clearpage

\clearpage

\clearpage



\begin{thebibliography}{9}
\addcontentsline{toc}{section}{\refname}




\bibitem{Beckner75} W.~Beckner, {Inequalities in Fourier analysis}, {\em Ann.~of~Math.} {\bf 102} (1975), 159--182.

\bibitem{BL} J.~Bergh, J.~L\"ofstr\"om, {\em Interpolation spaces. An introduction}.
Springer, 1976.

\bibitem{Bloom} S.~Bloom, 
Hardy integral estimates for the Laplace transform.
{\em Proc.\ Amer.\ Math.\ Soc.} {\bf 116} (1992), no.~2, 417--426.




\bibitem{Duren1970} P.L.~Duren, {\em Theory of $H^p$ spaces}. Academic Press, 1970.


\bibitem{Ed} R.E.~Edwards, {\em Fourier series. A modern introduction. Vol.~2}. Springer, 1982.

\bibitem{EvGa} L.C.~Evans, R.F.~Gariepy, {\em Measure theory and fine properties of funcions}, CRC Press, 2000. 

\bibitem{Garling} D.H.J.~Garling, {\em Inequalities}, Cambridge Univ.~Press, 2007.

\bibitem{GelOlm} B.~Gelbaum, J.~Olmsted, {\em Counterexamples in analysis}, Holden Day, 1964.

\bibitem{Grafakos}
L.~Grafakos, {\em Classical and Modern Fourier Analysis}. Pearson Education Inc., 2004.

\bibitem{Hardy33}
G.H.~Hardy, The constants of certain inequalities, 
{\em J.~London Math.~Soc}. {\bf 8} (1933), 114--119.


\bibitem{HLP}
G.H.~Hardy, J.E.~Littlewood, G.~Polya, {\em Inequalities}. Cambridge
Univ.\ Press, 1934.


\bibitem{Holmstedt} T.~Holmstedt, Interpolation of quasi-normed spaces. {\em Math.~Scand}. {\bf 26}
(1970), 177--199.

\bibitem{ALPDO}
L.~H\"ormander, {\em The analysis of linear partial differential operators. Vol.~1: Theory of distributions and Fourier analysis}. 
Springer, 1990. 



\bibitem{Merzon09}
A.E.~Merzon, F.-O.~Speck and T.J.~Villalba-Vega. On the weak
solution of the Neumann problem for the 2D Helmholtz equation in a convex cone and $H^s$ regularity. {\em Mathematical Methods in the Applied Sciences} {\bf 34} (2011), 24--43.

\bibitem{Oklander} E.T.~Oklander, 
$L_{pq}$ interpolators and the theorem of Marcinkiewicz,
{\em Bull.\ Amer.\ Math.\ Soc}. {\bf 72} (1966), 49--53.


\bibitem{RS} M.~Reed, B.~Simon, {\em Methods of modern mathematical physics Vol.~1: Functional analysis}.
Academic Press, 1972.


\bibitem{Sed2005} A.M.~Sedletskii, {\em Classes of analytic Fourier transforms, and exponential approximations}.
Fizmatlit, Moscow, 2005, in Russian. 

\bibitem{SM-CMA10} S.~Sadov, A.~Merzon,
$L^2$-estimates for the Laplace transform along a family of hyperbolas in the right half-plane,
in: Proceedings of ``Analysis, Mathematical Physics and Applications''
(Ixtapa, Mexico, March 1--5, 2010), 
{\em Comm.\ in Math.\ Analysis},
Conference 03 (2011), 204--208.

\bibitem{T}
E.~Titchmarsch, {\em Introduction to the theory of Fourier integrals},
Clarendon Press, 1948.

\bibitem{Sqv}
E.~Setterqvist. Unitary equivalence: a new approach to the Laplace operator and the Hardy operator.
{\em M.Sc.\ thesis}, Lule\o{a} Univ.\ of Technology, 2005. URL:
{\tt http://epubl.ltu.se/1402-1617/2005/329/LTU-EX-05329-SE.pdf} (accessed May 6, 2009).





\bibitem{ZM} P.~Zhevandrov, A.~Merzon, On the Neumann problem for
the Helmholtz equation in a plane angle.
{\em Mathematical Methods in  the Applied Sciences}. {\bf 23} (2000), 1401--1446.

\end{thebibliography}
\end{document}